\documentclass[11pt]{article} 
\usepackage[margin=1in]{geometry} 
\usepackage{setspace}
\usepackage[autostyle, english = american]{csquotes}

\usepackage[normalem]{ulem} 
\usepackage{hyperref}
\usepackage[english]{babel}

\usepackage{amsthm}
\usepackage{amsfonts}
\usepackage{newtxmath}
\usepackage{amsmath}

\newtheorem{proposition}{Proposition}
\newtheorem{theorem}{Theorem}

\newtheorem{corollary}{Corollary}
\newtheorem{claim}{Claim}

\usepackage{endnotes}
\let\footnote=\endnote

\usepackage{natbib}
\bibpunct[, ]{(}{)}{,}{a}{}{,}%
\usepackage{url}
\usepackage{nicefrac}
\usepackage{amsmath}
\usepackage{hyperref}
\usepackage{breakurl}
\newcommand{\LeftEqNo}{\let\veqno\@@leqno}
\usepackage{environ}
\usepackage{xcolor}
\usepackage{float}
\usepackage{graphicx}
\usepackage{amsmath}
\usepackage[noend]{algpseudocode}
\usepackage{algorithm}
\usepackage{tikz}
\usetikzlibrary{positioning,calc}
\usetikzlibrary{arrows,decorations.markings}
\usepackage{float}
\usepackage{url}
\usepackage{graphicx}
\usepackage{array}
\usepackage{enumerate}
\usepackage{caption}
\usepackage{subcaption}
\usepackage{comment}
\usepackage{lscape}
\usepackage{slashbox}

\usepackage[colorinlistoftodos,prependcaption,textsize=tiny, textwidth=\marginparwidth]{todonotes}

\usepackage[normalem]{ulem}

\usepackage{lscape}
\newcolumntype{P}[1]{>{\centering\arraybackslash}p{#1}}

\newenvironment{Proof}[1][Proof]{\begin{trivlist}
\item[\hskip \labelsep {\bfseries #1}]}{\end{trivlist}}

\newcommand{\set}[1]{\left\{#1\right\}}

\newcommand{\cuts}{\mathcal{C}} 
\newcommand{\N}{\mathcal{N}} 

\newcommand{\diff}{\delta}

\hypersetup{
colorlinks,
linkcolor={red!50!black},
citecolor={blue!50!black},
urlcolor={blue!80!black}
}
\makeatletter
\NewEnviron{Lalign}{\tagsleft@true\begin{align}\BODY\end{align}}
\makeatother

\newcolumntype{C}[1]{>{\centering\let\newline\\\arraybackslash\hspace{0pt}}m{#1}}

\begin{document}

\title{Optimizing the expected maximum of two linear functions defined on a multivariate Gaussian distribution}

\author{David Bergman \\ \small Department of Operations and Information Management, School of Business, University of Connecticut
\and
Carlos Cardonha \\ \small Department of Operations and Information Management, School of Business, University of Connecticut
\and 
Jason Imbrogno \\ \small Department of Finance, Economics, and Data Analytics,  
University of North Alabama
\and
Leonardo Lozano  \\ \small Operations, Business Analytics \& Information Systems, University of Cincinnati
}
\date{}
\maketitle

\abstract{%
We study stochastic optimization problems with objective function given by the expectation of the maximum of two linear functions defined on the component random variables of a multivariate Gaussian distribution. We consider random variables that are arbitrarily correlated, and we show that the problem is NP-hard even if the space of feasible solutions is unconstrained. We exploit a closed-form expression for the objective function from the literature to construct a cutting-plane algorithm that can be seen as an extension of the integer L-shaped method for a highly nonlinear function, which includes the evaluation of the c.d.f and p.d.f of a standard normal random variable with decision variables as part of the arguments. To exhibit the model's applicability, we consider two featured applications. The first is daily fantasy sports, where the algorithm identifies entries with positive returns during the 2018-2019 National Football League season. The second is a special case of makespan minimization for two parallel machines and jobs with uncertain processing times; for the special case where the jobs are uncorrelated, we prove the equivalence between its deterministic and stochastic versions and show that our algorithm can deliver a constant-factor approximation guarantee for the problem. The results of our computational evaluation involving synthetic and real-world data suggest that our discretization and upper bounding techniques lead to significant computational improvements and that the proposed algorithm outperforms sub-optimal solutions approaches. 
}



\maketitle
%




\section{Introduction}
\label{sec:introduction}

Consider two zero-mean jointly normal random variables, and the expectation of their maximum.  Should they be uncorrelated, the expectation of the maximum is $\frac{1}{\sqrt{\pi}}$.  With correlation $\frac{1}{2}$ the expectation is $\frac{1}{\sqrt{2\pi}}$, and with correlation $-\frac{1}{2}$ the expectation is $\frac{\sqrt{3}}{\sqrt{2\pi}}$.  For arbitrary jointly normal random variables $Z_1,Z_2$ the maximum is  as follows (see~\cite{nadarajah2008exact}):
\begin{eqnarray}
 \mathbb{E}\left[\max \{Z_1,Z_2\}\right] =&
\mathbb{E} 
\left[ 
	Z_1 
\right]
\Phi
\left(
	\frac
    	{     
\mathbb{E}[Z_1] - \mathbb{E}[Z_2]
        }
    	{\sqrt{\sigma^2(Z_1) + \sigma^2(Z_2) - 2 \mathrm{cov}(Z_1, Z_2)}}
\right)
+  \nonumber \\
& \mathbb{E} 
\left[ 
	Z_2
\right]
\Phi
\left(
	\frac
    	{     
\mathbb{E}[Z_2] - \mathbb{E}[Z_1]
        }
    	{\sqrt{\sigma^2(Z_1) + \sigma^2(Z_2) - 2 \mathrm{cov}(Z_1, Z_2)}}
\right) + \nonumber \\
& \sqrt{\sigma^2(Z_1) + \sigma^2(Z_2) - 2\mathrm{cov}(Z_1, Z_2)}
\phi
\left(
	\frac
    	{  
             E[Z_1] - E[Z_2]
        }
    	{\sqrt{\sigma^2(Z_1) + \sigma^2(Z_2) - 2 \mathrm{cov}(Z_1, Z_2)}}
\right),& \nonumber
\end{eqnarray}
where $\Phi$ and $\phi$ are the c.d.f. and p.d.f. of a standard normal random variable, respectively. The evaluation of the expectation requires the calculation of the highly nonlinear c.d.f. and p.d.f. with arguments that are composed of ratios of linear expressions over square roots of linear expressions (and even the product of these expressions). Consider now an optimization problem, where the objective function consists of the expression above, 
but the selection of $Z_1$ and $Z_2$ are subject to decisions in the problem.  This paper studies this decision making setting, describes real-world problems that can be modeled within this context, and explores an efficient algorithm for solving the resulting optimization problems. More formally, we investigate a class of stochastic optimization problems defined over a multivariate Gaussian distribution.  In particular, we study the problem of optimizing the expected value of the maximum of two linear functions defined on the component random variables of a multivariate Gaussian distribution, allowing for correlation among the component random variables. 

Real-world problems that can be modelled within this context abound, especially if we consider the  connection to order statistics, which is used to model types of auctions  where the final prices are determined by the first- or second-highest bidder \citep{brown1986using}; 
insurance premium determination, where companies use order statistics  to  determine policies for joint-life insurance \citep{
city20705}; failure models for wireless communication networks~\citep{yang_alouini_2011}; and risk management \citep{Koutras2018}.
We study two applications in this paper\textemdash  (1) expected score maximization in Daily Fantasy Sports (DFS) and (2) makespan minimization for 2 parallel machines with stochastic processing times. These two applications highlight how both objectives (maximizing and minimizing) are important in different contexts. 

 Exact algorithms for optimizing expressions composed of the maximum of two different stochastic functions
 is challenging due to the complexity of evaluating or even estimating the expectation of order statistics \citep{david2004order,Bertsimas:2006:TBE:1187913.1187922, evans2006distribution}.  Specific to the problem class studied in this paper, closed-form expressions for the expected value of the maximum and the minimum are known (see~\cite{nadarajah2008exact}). 
We show how one can formulate the problem as a binary optimization problem and  investigate exact computational approaches for solving the model. Our approach consists of 
an extension of the integer L-shaped method 
combined with linearization and discretization techniques tailored for the problem. 

Our contributions are the following:
\begin{itemize}
    \item We show that the underlying optimization problem is NP-hard even in scenarios where the set of solutions is unconstrained;
    
    \item We develop an exact optimization algorithm that allows for the incorporation of linear constraints and establish theoretical results on the quality of the upper bounds obtained;

    \item We study two featured applications, aimed at exhibiting the generality of the class of problems and to showcase the scalability of the proposed algorithms. 
    
    \item  For~$P2|p \sim \N(\mu,\Sigma),\rho_{j,j'}=0|\mathbb{E}[C_{max}]$, a machine scheduling problem investigated in this article, we show that a single iteration of our algorithm can deliver a solution with constant-factor approximation guarantees, 
    and we also prove that this problem 
    is equivalent to its deterministic counterpart (where processing times are replaced by the mean of the respective distributions); and
    
    
    \item We conduct a comprehensive computational study both on synthetic problem instances and real-world data in order to show the computational performance and the quality of the results delivered by our algorithms.
  
\end{itemize}
    
The paper is organized as follows. 
Section~\ref{sec:probDesc} formally defines the class of problems we study. Section~\ref{sec:complexity} presents computational complexity results. In Section~\ref{sec:Exact} we present an exact cutting-plane algorithmic framework. Numerical experiments on synthetic instances are presented and discussed in Section \ref{sec:compExperiments}. Sections  \ref{sec:DFF} and \ref{sec:scheduling} present our two featured applications. Section \ref{sec:conclusion} concludes the paper and discusses future work.

\section{Problem Description and Featured Applications}
 \label{sec:probDesc}

We consider a collection of random variables $\pmb{Y} = \left( Y_1, \ldots, Y_n \right)^T$ following a multivariate Gaussian distribution $\pmb{Y} \sim \mathcal{N}(\pmb{\mu},\pmb{\Sigma})$. The random variables $Y_1, \ldots, Y_n$ will be referred to as \emph{component} random variables which have component-wise expectations $\mu_1, \ldots, \mu_n$.  The $n \times n$ covariance matrix $\pmb{\Sigma}$ has elements $\pmb{\Sigma}_{i,j} = \textnormal{Cov}(Y_i, Y_j)$.
 
We study the class of problems
\begin{equation}
\tag{P}
\label{eqn:optProblem}
\max\limits_{x \in \Omega \subseteq \{0,1\}^{2 \times n}}   \textnormal{(or $\min$)} \quad \mathbb{E} \left[ \max \{ Z_1(x),Z_2(x) \} \right], 
\end{equation}
where for $i=1,2$, $Z_i(x) = \sum_{j=1}^n Y_j x_{i,j} $ with $x_{i,j} \in \set{0,1}$ for $j = 1, \ldots, n$.  The set $\Omega$ defines the feasible region for the decision variables, assumed in this paper to be composed solely of  linear constraints. 
In order to exhibit the expansive generalizability of this class of problems, we present two featured applications.

\noindent \textbf{Featured Application 1: Daily Fantasy Sports}

In certain DFS competitions, participants select up to~2 fantasy entries, with each entry being composed of a set of players who will participate in upcoming sporting events. The selection of players is subject to roster and budget constraints; more details about these and other constraints are presented in Section~\ref{sec:DFF}. Each fantasy entry receives points based on the actual performance of the players selected in the sporting events. Entries are ranked according to the total fantasy points scored, and the payout of each entry depends on its position in this rank.  The payout structures are top-heavy, with a small fraction of the entries receiving substantial amounts, including the top-scoring entry receiving approximately 25\% of the total entry fees paid.
 
Suppose there are $n$ players that can be selected for any of a participant's $2$ entries.  Each player will score a random number of fantasy points $Y_j$.  Letting $x_{i,j}$ indicate the selection of player $j$ for entry $i$, for $i = 1, 2$ and $j = 1, \ldots, n$, one can model the selection problem for a participant as \begin{align}
& \max && \mathbb{E} \left[ \max \left\{ Z_1(x), Z_2(x) \right\}  \right]  \label{fa:dff} \tag{DFS} \\
& \textnormal{s.t.} && Z_i(x) = \sum_{j=1}^n Y_j x_{i,j} && i = 1, 2 \nonumber \\
&&& x \in \Omega &&\nonumber 
\end{align}
The correlation between the points received for a pair of players can be significant, for example between a quarterback and a wide receiver in football, since most fantasy points that a wide receiver receives will generally be associated with fantasy points for the quarterback on the same team.  The constraints in $\Omega$ define conditions which enforce what configurations of players make a legal entry, which will be explained in more detail later and can be modeled through linear constraints.

Algorithmic sports betting recently became a topic of interest in the operations management and operations research literature~\citep{KapGar2001,ClaLet07,haugh2021play}. Some examples involve selecting multiple entries for maximizing the expected score of the maximum scoring entry in both National Football League (NFL) survival pools (\cite{bergman2017surviving}) and DFS (\cite{hunter2016picking,haugh2021play}).

\noindent \textbf{Featured Application 2: Makespan Minimization  with Stochastic Processing Times}

Consider a set of $n$ jobs with processing times~$Y_j$ drawn from a multivariate Gaussian distribution to be partitioned for execution on~$2$ parallel (identitical) machines such that the makespan (the time at which the last job completes) is minimized; this problems is represented by~$P2|p \sim \N(\mu,\Sigma)|\mathbb{E}[C_{max}]$
in the notation of~\cite{graham1979optimization}. 
Letting binary variable $x_{i,j}$ indicate if job $j$ is assigned to machine $i$ for $i=1, 2$ and $j=1, \ldots, n$, the problem can be formulated as
\begin{align}
& \min && \mathbb{E} \left[ \max \left\{ Z_1(x),Z_2(x) \right\} \right]  \label{fa:ms} \tag{MS} \\
& \textnormal{s.t.} && Z_i(x) = \sum_{j=1}^n Y_j x_{i,j} && i = 1, 2 \nonumber \\
&&& x_{1,j} + x_{2,j} = 1 && j=1, \ldots, n \nonumber \\
&&& x_{i,j} \in \set{0,1} && i=1, 2, \ j = 1, \ldots, n. \nonumber
\end{align}
The objective seeks to minimize the maximum expected completion time of all machines. Note that no assumption is presumed on the correlation between the processing time of the jobs. We incorporate $\rho_{j,j'}=0$ to the notation  to represent scenarios where processing times are uncorrelated.

There is a vast literature on stochastic scheduling~(see~\cite{nino2009stochastic}). In particular, \cite{coffman1987minimizing} characterize optimal solutions for makespan minimization in scenarios involving 2 or 3 machines and jobs with exponentially distributed processing times. \cite{pinedo2005planning}  
shows the connection between a stochastic flow shop scheduling problem and a deterministic traveling salesman problem; we present a similar result involving~$P2|p \sim \N(\mu,\Sigma),\rho_{j,j'}=0|\mathbb{E}[C_{max}]$ in this~work.

\section{Computational Complexity}
\label{sec:complexity}
Consider random variables $\pmb{Y} = \left( Y_1, \ldots, Y_n \right)^T$ following a multivariate Gaussian distribution $\pmb{Y} \sim \mathcal{N}(\pmb{\mu},\pmb{\Sigma})$ and note that  all component random variables~$Y_{j}$ are normally distributed and potentially correlated, i.e., 
$Y_j \sim \N(\mu_{j},\sigma^2_j)$,
with mean $\mu_j$ and variance $\sigma^2_j = \pmb{\Sigma}_{j,j}$, for all $j \in \set{1, \dots, n}$.

\begin{theorem}\label{thm:hard}
Optimization problem \ref{eqn:optProblem} is NP-hard even if $\Omega = \{0,1\}^{2 \times n}$.
\end{theorem}

Before proceeding with the proof, we recall some known results that will be relevant in this section and throughout the manuscript. Since $\pmb{Y}$ follows a multivariate Gaussian distribution, $Z_i(x) = \sum\limits_{j=1}^n Y_j x_{i,j} $ is normally distributed for all $i=1,2$. Moreover, for any $x \in \{0,1\}^{2 \times n}$, the mean and variance of~$Z_{1}(x)$ and~$Z_{2}(x)$, as well as their covariances (and correlations), are given~by:
\begin{align}
& \mathbb{E} \left[Z_i(x) \right] = \sum_{j = 1}^n \mu_j  x_{i,j} && \forall i \in \{1,2\} \label{expectedVal} \\
& \sigma^2 \left(Z_i(x) \right) =
\sum_{j=1}^n \sigma^2_j  x_{i,j}
+ 
2 \sum_{1 \leq j < j' \leq n} \mathrm{cov}(Y_{j} , Y_{j'} )  x_{i,j}  x_{i,j'} 
&& 
\forall i \in \{1,2\} \label{variance} 
\\
&
\mathrm{cov} \left( Z_{1} (x) ,  Z_{2} (x) \right) = 
\sum_{j = 1}^{n} \sum_{j' = 1}^n  \mathrm{cov} \left( Y_{j} , Y_{j'}  \right) x_{1,j} x_{2,j'}.
&&
\label{covariance}
\end{align}

Let $\phi(\cdot)$ and $\Phi(\cdot)$  be the probability density function (p.d.f.) and the cumulative distribution function (c.d.f.), respectively, of standard normal random variables; i.e., for $w \in (- \infty , \infty )$,
\begin{align*}
& \phi (w) = \frac{1}{\sqrt[]{2 \pi}} e ^ {\frac{-w^2}{2}},&
& \Phi (w) = \int_{- \infty}^w \phi (u) du.
\end{align*}

An exact expression is known for $\mathbb{E} \left[ \max \{ Z_1(x),Z_2(x) \} \right]$ (see~\cite{nadarajah2008exact} and~\cite{clark1961greatest}):
\begin{equation}\label{eq:objective}
\mathbb{E} 
\left[ 
	Z_1 (x)
\right]
\Phi
\left(
	\frac
    	{     
\delta(x)
        }
    	{\theta(x)}
\right)
+ 
\mathbb{E} 
\left[ 
	Z_2 (x)
\right]
\Phi
\left(
	\frac
    	{    
-\delta(x)
        }
    	{\theta(x)}
\right) +
\theta(x)
\phi
\left(
	\frac
    	{  
             \delta(x)
        }
    	{\theta(x)}
\right),
\end{equation}
where
\begin{equation*}
    \delta(x) =         	
    \mathbb{E} 
          		\left[ 
              	Z_1 (x)
          	\right]
            - 
        	\mathbb{E} 
          		\left[ 
              	Z_2 (x)
          	\right],
\end{equation*}
and 
\begin{equation}
\theta(x)\label{eq:theta}
=
\sqrt
{
	\sigma^2(Z_1(x))
    +
	\sigma^2(Z_2(x))
    -
    2 \mathrm{cov} \left(Z_1(x) , Z_2(x)\right)
}.
\end{equation}
In order to simplify notation, we assume without loss of generality that $\mathbb{E}\left[Z_1 (x)\right] \geq \mathbb{E} \left[Z_2 (x)\right]$.

As we can see from Expression~\ref{eq:objective}, due to the (possible) dependency among the random variables $Z_1(x)$ and $Z_2(x)$, the expectation of their maximum is a highly nonlinear function. Note that if $\theta(x)$ is 0, the ratio in the argument of the c.d.f. may be undefined, but the expression holds true if we assume $\Phi\left( \frac{a}{0} \right) $ to be 0 if $a < 0$, $\frac{1}{2}$ if $a = 0$, and 1 if $a > 0$.

Equipped with this expression, we proceed with the proof of Theorem~\ref{thm:hard}.

\begin{proof}{Proof of Theorem~\ref{thm:hard}.}
The result follows from a reduction from the minimum weighted cut problem, which  is known to be NP-hard (\cite{mccormick2003easy}).   Let $G = (V,E)$ be an undirected graph with integer weights $w_e$ on each edge $e \in E$. An $(S,T)$-\emph{cut} of $G$ is a 2-partition of $V$, and its  \emph{weight}~$w(S,T)$ is the sum of the weights of the edges ``crossing'' the cut, i.e., $w(S,T) = \sum\limits_{e : e \cap S, e \cap T \neq \emptyset} w_e$.  In the decision version of the  problem, we are also given a constant $K$  and wish to know if there exists a~$(S,T)$-cut  of~$G$ such that $w(S,T) \leq K$.  

We create an instance of \ref{eqn:optProblem} as follows.  Assuming there are $n$ vertices in $G$, every vertex $j \in V$ is associated with a normally distributed random variable~$Y_j \sim \N(0,1)$, for $j = 1, \ldots, n$. The covariance between $Y_{j}$ and $Y_{j'}$ is $\mathrm{cov}(Y_{j},Y_{j'}) = \frac{w_{\{j,j'\}}}{4 M+1}$, where $M = \sum\limits_{e \in E} |w_e|$. As a result, we have that $\sum\limits_{j = 1}^n \sum\limits_{\substack{j' =  1\\j' \neq j}}^n \mathrm{cov}(Y_{j},Y_{j'}) \leq \frac{2M}{4M+1} < \frac{1}{2}$ and $|\sigma^2_{j}| \geq \sum\limits_{j' \neq j}|\mathrm{cov}(Y_{j},Y_{j'})|$ for all~$j= 1, \ldots, n$. One can then construct a symmetric and diagonally dominant (consequently, positive semi-definite) matrix~$\Sigma$ whose columns and rows are indexed by variables~$Y_j$. It follows that~$\Sigma$ is a valid covariance matrix.
Finally, let $\Omega = \{0,1\}^{2 \times n}$, i.e., the set of feasible solutions is unconstrained.

By construction, $\mathbb{E}[Z_1(x)] = \mathbb{E}[Z_2(x)] = 0$, so the expression for $\mathbb{E}[Z_{(2)}(x)]$ reduces to~$\theta(x) \frac{1}{\sqrt{2 \pi }}$.
As~$n \geq 1$ and all variances are equal to 1, it follows that at optimality $\theta(x) \geq 1$ (this value is achieved if we assign one item to the first set and leave the second set empty), so any~$x$ that maximizes $\theta(x)^2$ also maximizes $\theta(x)$. Therefore, our problem is  equivalent to
\[
 \max_{x \in \Omega} \ \theta(x)^2 =  \sigma^2(Z_1(x)) + \sigma^2(Z_2(x)) - 2 \mathrm{cov} (Z_1(x),Z_2(x)).
\]
By expanding the terms of the last expression and replacing all variances for 1, $\theta(x)^2$ becomes  
\begin{equation}
\begin{aligned}\label{thm:theta2}
&\sum_{j = 1}^n (x_{1,j} + x_{2,j} - 2x_{1,j}x_{2,j} ) \, +  \\
 & \quad \quad  2
\left[
\sum_{j = 1}^{n-1} \sum_{j' = j + 1}^n \mathrm{cov}(Y_{j}, Y_{j'}) (x_{1,j}  x_{1,j'} 
+ x_{2,j}  x_{2,j'})
- \sum_{j = 1}^{n} \sum_{\substack{j' =  1\\j' \neq j}}^n \mathrm{cov}(Y_{j}, Y_{j'})x_{1,j}  x_{2,j'} 
\right]
\end{aligned}
\end{equation}

\begin{claim}\label{claim1}
In any optimal solution for $\max_{x \in \Omega} \theta(x)^2 $, $x_{1,j} + x_{2,j} = 1$ for $j = 1, \ldots, n$.
\end{claim} 
\begin{Proof}
Let~$A(x)$ denote the sum within the brackets in~\ref{thm:theta2}; 
$A(x)$ belongs to the interval defined by
 $   \pm   \sum\limits_{j = 1}^n \sum_{\substack{j' =  1,j' \neq j}}^n \mathrm{cov}(Y_{j},Y_{j'})$,
because each covariance term appears at most twice with a positive coefficient and at most twice with a negative coefficient.
By construction, $ \sum\limits_{j = 1}^n \sum_{\substack{j' =  1\\j' \neq j}}^n \mathrm{cov}(Y_{j},Y_{j'}) < \frac{1}{2}$, so
we have that  $-1 < 2A(x) < 1 $.  Therefore, any optimal solution to $\max_{x \in \Omega} \ \theta(x)^2 $ also optimizes
\[
\max_{x \in \Omega} \left\{ \sum_{j = 1}^n x_{1,j}
+  \sum_{j = 1}^n x_{2,j}
-
2 \sum_{j=1}^n x_{1,j} x_{2,j} \right\},
\]
since the absolute value of each of the coefficients in this expression is greater than or equal to 1.   
Finally, note that this expression is maximized if and only if $x_{1,j} + x_{2,j} = 1$, as desired.
$\blacksquare$
\end{Proof}
It follows from  Claim~\ref{claim1} that 
$\sum\limits_{j = 1}^n 
(
x_{1,j}
+ x_{2,j}
-
2 x_{1,j} x_{2,j} )
= n$, so the problem reduces to
\begin{align}
\label{eq:rewrite1}
 \max_{x \in \Omega}\quad  &  
\sum\limits_{j = 1}^{n-1} \sum\limits_{j' = j + 1}^n \mathrm{cov}(Y_{j}, Y_{j'}) (x_{1,j}  x_{1,j'} 
+ x_{2,j}  x_{2,j'})
- \sum\limits_{j = 1}^{n} \sum\limits_{\substack{j' =  1\\j' \neq j}}^n \mathrm{cov}(Y_{j}, Y_{j'})x_{1,j}  x_{2,j'} &\\
 \textnormal{s.t.}\quad &  x_{1,j} + x_{2,j} = 1 & j \in [n] \\
&  x_{i,j} \in \{0,1\} & i \in [2],j \in [n]. 
\end{align}
Consider the following optimization problem:
\begin{equation}
\label{eq:rewrite2}
\begin{aligned}
& \min    && h(x) = \sum_{j = 1}^{n} \sum\limits_{\substack{j' =  1,j' \neq j}}^n \mathrm{cov}(Y_{j}, Y_{j'})x_{1,j} x_{2,j'}  \\
& \textnormal{s.t.} &&
x_{1,j} + x_{2,j} = 1 && j \in [n] \\
&&& x_{i,j} \in \{0,1\} && i \in [2], \ j \in [n]. 
\end{aligned}
\end{equation}

\begin{claim}
An optimal solution to optimization problem~(\ref{eq:rewrite2}) is also optimal to
problem~(\ref{eq:rewrite1}).
\end{claim}
\begin{Proof}
Both optimization problems have the same set of feasible solutions~$\Omega'$.  Let $x'$ and $x''$ be two feasible solutions with $ h(x') < h(x'') $.  Showing that $\theta(x')^2 > \theta(x'')^2$ establishes the claim.  To show this, we first note that for any feasible solution $\tilde{x}$, 
\begin{equation}
\label{eq:constant}
\begin{aligned}
&  \sum_{j = 1}^{n-1} \sum_{j' = j + 1}^n \mathrm{cov}(Y_{j}, Y_{j'})\tilde{x}_{1,j}  \tilde{x}_{1,j'}
+ \sum_{j = 1}^{n-1} \sum_{j' = j + 1}^n    \mathrm{cov}(Y_{j}, Y_{j'})\tilde{x}_{2,j}  \tilde{x}_{2,j'} \\
& + \sum_{j = 1}^{n} \sum_{j' \neq j
}^p \mathrm{cov}(Y_{j}, Y_{j'})\tilde{x}_{1,j}  \tilde{x}_{2,j'} = \sum_{j = 1}^{n-1} \sum_{j'=j + 1}^{n} \mathrm{cov}(Y_{j}, Y_{j'}).
\end{aligned}
\end{equation}
This follows because for any two indices $j \neq j'$, the covariance term $\mathrm{cov}(Y_{j},Y_{j'})$ is counted in exactly one of the three terms in the left-hand size of equation~(\ref{eq:constant}):
\begin{enumerate}
    \item If $\tilde{x}_{1,j} = \tilde{x}_{1,j'} = 1 \rightarrow$ $\mathrm{cov}(Y_{j}, Y_{j'})$ is counted only in the first term; 
    \item If $\tilde{x}_{2,j} = \tilde{x}_{2,j'} = 1 \rightarrow$ $\mathrm{cov}(Y_{j}, Y_{j'})$ is counted only in the second term;   
    \item If $\tilde{x}_{1,j} = 1, \tilde{x}_{2,j'} = 1 \rightarrow$ $\mathrm{cov}(Y_{j}, Y_{j'})$ is counted only in the third term; and 
    \item If $\tilde{x}_{2,j} = 1, \tilde{x}_{1,j'} = 1 \rightarrow$ $\mathrm{cov}(Y_{j}, Y_{j'})$ is counted only in the third term.
\end{enumerate}
Finally, because $\tilde{x}_{1,j} + \tilde{x}_{2,j} = 1$ and $\tilde{x}_{1,j'} + \tilde{x}_{2,j'} = 1$, it follows that the list above is exhaustive and contains all possible assignments of~$j$ and~$j'$. Therefore, 
\begin{eqnarray*}
h(x')  < h(x'') &\implies& 
-h(x')  > -h(x'') \implies  \\
\sum_{j = 1}^{n-1} \sum_{j'=j+1}^{n} \mathrm{cov}(Y_{j}, Y_{j'}) - h(x')
&>&
\sum_{j = 1}^{n-1} \sum_{j'=j+1}^{n} \mathrm{cov}(Y_{j}, Y_{j'}) - h(x''),
\end{eqnarray*}
which implies that $\theta(x')^2 > \theta(x'')^2$, 
as desired. 
$\blacksquare$
\end{Proof}
We conclude that an optimal solution to~(\ref{eq:rewrite2}) is also optimal to the original problem.
We now show a one-to-one mapping between solutions of~(\ref{eq:rewrite2}) and $(S,T)$-cuts of $G$.  Namely, for a feasible solution $x'$, we associate it with the $(S(x'),T(x'))$-cut
defined by $x'_{1,j} = 1 \leftrightarrow j \in S(x')$ and $x'_{2,j} = 1 \leftrightarrow j \in T(x')$. Additionally, $w(S(x'),T(x')) \leq K$ if and only if $h(x') \leq \frac{K}{4M + 1}$. This follows because 
\begin{equation*}
\begin{aligned}
w(S(x'),T(x')) & = \sum_{j \in S(x')} \sum_{j' \in T(x')} w_{\{j,j'\} }  = 
\sum_{j \in S(x')} \sum_{j' \in T(x')}
\left( 4M + 1 \right) 
\mathrm{cov}\left(Y_{j}, Y_{j'}\right) \\
& = 
\left( 4M + 1 \right) 
\sum_{j = 1}^n \sum_{j' = 1, j' \neq j}^n 
\mathrm{cov}\left(Y_{j}, Y_{j'}\right)  = 
\left( 4M + 1 \right) h(x').
\end{aligned}
\end{equation*}
It follows that~(\ref{eq:rewrite2}) is equivalent to the  minimum weighted cut problem, as desired.
$\blacksquare$
\end{proof}


\section{Exact Solution Algorithm} \label{sec:Exact}
We propose an exact solution algorithm for problem \ref{eqn:optProblem} and focus on the maximization of the objective function. We discuss in the Appendix how to adapt the algorithm for the minimization case. Omitted proofs can be found in the Appendix.

\subsection{A Cutting-Plane Framework}
\label{sec:Algorithm}
The objective function presented in \eqref{eq:objective} is highly nonlinear, thus making the application of direct formulations combined with off-the-shelf solvers unlikely to be successful. We present an exact cutting-plane algorithm based on a mixed-integer linear program (MILP) referred to as the \emph{relaxed master problem}. This MILP provides upper bounds on the optimal objective value and is iteratively updated through the inclusion of cuts. Lower bounds are obtained through the evaluation of~$\mathbb{E} \left[ \max \{ Z_1(x),Z_2(x) \} \right]$ on the feasible solutions obtained from the optimization of the relaxed master problem. The algorithm stops once it finds a provable optimal solution.

Our approach to solve the problem is presented in Algorithm~\ref{a1}. A key component of our algorithm is the construction of an upper-bounding function for~$\mathbb{E} \left[ \max \{ Z_1(x),Z_2(x) \} \right]$ defined over the set~$\Omega$ of feasible solutions. Namely, we wish to work with a function~$g(x)$ such that 
$
g(x) \geq \mathbb{E} \left[ \max \{ Z_1(x),Z_2(x) \} \right]$ for each~$x$ in~$\Omega$.
Given~$g(x)$, the relaxed master problem can be stated as
\begin{equation}
\tag{RMP}
\label{eqn:RMP}
\bar{z} = \max\limits_{x \in \Omega} g(x),
\end{equation}
and $\bar{z}$ provides an upper bound on the optimal objective value of problem~\ref{eqn:optProblem}. 

Algorithm~\ref{a1} solves problem RMP iteratively, adding \emph{no-good constraints} (or simply \emph{cuts}) to a set~$\cuts$, which are incorporated to RMP in order to prune previously explored solutions~\citep{balas72}. RMP$(\cuts)$ denotes RMP further constrained by~$\cuts$, so that a solution~$x$ for RMP$(\cuts)$ belongs to $\Omega$ and  satisfies all no-good constraints of~$\cuts$. In each iteration, after finding an optimal  solution~$\hat{x}$ for RMP$(\cuts)$, Algorithm~\ref{a1}  adds the following cut~$c(\hat{x})$ to~$\cuts$. 
\begin{align}\label{noGood}
    \sum_{i\in\set{1,2}}\sum_{j\in \set{1,\dots,n}:\hat{x}_{i,j}=1} x_{i,j} -\sum_{i\in\set{1,2}}\sum_{j\in \set{1,\dots,n}:\hat{x}_{i,j}=0} x_{i,j} \leq  \sum_{i\in\set{1,2}}\sum_{j\in \set{1,\dots,n}} \hat{x}_{i,j}-1.
\end{align}
The only solution in~$\Omega$ that violates~$c(\hat{x})$ is~$\hat{x}$. Therefore, every solution in~$\Omega$ is explored at most once, and as~$\Omega$ is finite, Algorithm~\ref{a1} terminates after a finite number of steps.

Our cutting-plane algorithm keeps a lower bound~\textit{LB} and an upper bound~\textit{UB} for~$\mathbb{E} \left[ \max \{ Z_1(x),Z_2(x) \} \right]$, 
which are iteratively updated based on solutions~$\hat{x}$ of RMP($\cuts$). Upper bounds are given by~$g(\hat{x})$; these bounds  are non-decreasing, as any solution~$x'$ such that $g(x') > g(\hat{x})$ must have been explored in a previous iteration of Algorithm~\ref{a1}. Similarly, lower bounds are obtained through the exact evaluation of $\mathbb{E} \left[ \max \{ Z_1(\hat{x}),Z_2(\hat{x}) \} \right]$ (i.e., using Equation~\eqref{eq:objective}). Observe that $\mathbb{E} \left[ \max \{ Z_1(\hat{x}),Z_2(\hat{x}) \} \right]$ may be smaller than~$\mathbb{E} \left[ \max \{ Z_1(x'),Z_2(x') \} \right]$ for some previously explored solution~$x'$, so Algorithm~\ref{a1} needs to store the largest~$LB$ found in previous iterations. Algorithm~\ref{a1} terminates when $\textit{LB}$ becomes equal to~$\textit{UB}$ or when RMP($\cuts$) becomes infeasible. 

Algorithm \ref{a1} can be seen as an extension of the integer L-shaped method \citep{Laporte93}, dealing with the added difficulty of a highly nonlinear objective function. Similar cutting-plane algorithms have been extensively used in the context of two-stage stochastic programming with integer recourse \citep{
sen2006decomposition,angulo2016} and are closely related to the logic-based Benders' decomposition algorithm \citep{hooker2003logic}. As a result, our main contribution lies in the proposed linear upper bounding function described in Section \ref{sec:Bounds}.

\begin{algorithm}[H]
\begin{algorithmic}[1]

\State Set \textit{LB}$ = -\infty$, \textit{UB}$ = \infty$, $\mathcal{C} = \emptyset$, and incumbent solution $\bar{x} = 0$.

\State Optimize RMP$(\mathcal{C})$ to obtain~$\hat{x}$;  if the problem is infeasible, go to Step 6.  

\State 
Set \textit{UB} = $g(\hat{x})$.

\State 
If $\mathbb{E} \left[ \max \{ Z_1(\hat{x}),Z_2(\hat{x}) \} \right]$ $>$ \textit{LB}, set \textit{LB} = $\mathbb{E} \left[ \max \{ Z_1(\hat{x}),Z_2(\hat{x}) \} \right]$ and update incumbent $\bar{x} = \hat{x}$. 

\State If \textit{LB} = \textit{UB}, go to Step 6.  Otherwise, set~$\mathcal{C} = \mathcal{C} \cup \{c(\hat{x})\}$ and 
return to Step 2.
 
\State If \textit{LB} = $-\infty$, 
original problem is infeasible. Otherwise,
terminate with optimal solution 
$\bar{x}$.
\end{algorithmic}
\caption{A Cutting-Plane Algorithm}
\label{a1}
\end{algorithm}

\subsection{Baseline Approach for Obtaining Upper Bounds on the Objective Function}
\label{sec:Baseline}
The performance of our cutting-plane algorithm is directly tied to the quality of the bounds obtained from~$g(x)$ 
and the difficulty of solving the RMP. A bounding function that delivers accurate overestimations of $\mathbb{E} \left[ \max \{ Z_1(x),Z_2(x) \} \right]$ but requires an impractical amount of time to solve the resulting RMP is likely to result in poor performance of the cutting-plane algorithm, as the RMP is solved at every iteration of the algorithm. On the other hand, a bounding function for which the resulting RMP is easily solved but provides poor quality upper bounds is likely to result in a large number of iterations before closing the optimality gap.

We remark that $\mathbb{E} \left[ \max \{ Z_1(x),Z_2(x) \} \right] \not=  \max \{\mathbb{E} \left[Z_1(x)\right] ,\mathbb{E} \left[Z_2(x)\right] \}$, and in our setting it holds that $\mathbb{E} \left[ \max \{ Z_1(x),Z_2(x) \} \right] \geq  \max \{\mathbb{E} \left[Z_1(x)\right] ,\mathbb{E} \left[Z_2(x)\right] \}$. As a result, a simple bounding function is given by the following linear expression:
\begin{equation} \label{eq:simple}
    \max \{\mathbb{E} \left[Z_1(x)\right] ,\mathbb{E} \left[Z_2(x)\right] \} + M,
\end{equation}
where $M$ is a sufficiently large constant such that 
\begin{equation*} 
        M \geq \theta(x) \phi \left(\frac{\diff(x)}{\theta(x)} \right).
\end{equation*}
We formulate the corresponding RMP as a linear mixed-integer program (MIP) that can be solved with any off-the-shelf optimization software. Unfortunately, such a simple function yields poor quality bounds and virtually requires a 
complete enumeration of the solution space. In preliminary computational experiments, the cutting-plane algorithm is not able to solve to optimality a single problem instance from our test bed using \eqref{eq:simple} as the bounding function, resulting in high optimality gaps at the end of the time limit. As a result, we investigate more complex linear bounding functions, aiming to achieve a balance between the difficulty of solving the RMP (which we cast as a linear MIP) and the quality of the upper bounds obtained.    

A challenging task involved in obtaining bounds on the optimal objective is the evaluation of~$\theta(x)$, a nonlinear expression that appears in all terms of \eqref{eq:objective}, including the denominators of the c.d.f. and the p.d.f. of the standard normal distribution. To avoid the technical issues involved in the evaluation of~$\theta(x)$, we propose a baseline approach that evaluates~$\theta(x)^2$ exactly and defines an upper bounding function for~$\mathbb{E} \left[ \max \{ Z_1(x),Z_2(x) \} \right]$, based on a discretization of $\theta(x)^2$, where the value of $\theta(x)^2$ is 
\begin{equation}
\label{thetaSquared}
\theta(x)^2 = \sigma^2(Z_1(x)) + \sigma^2(Z_2(x))
    - 2 \mathrm{cov} \left(Z_1(x) , Z_2(x)\right).
\end{equation}
The exact value of~$\theta(x)^2$ is computed  via a McCormick linearization technique \citep{McCormick1976}. Let function $u_\theta(x)$ be such that $0  \leq \theta(x) \leq u_\theta(x)$. Our baseline upper-bounding function is given in Proposition \ref{Pbase}. 

\begin{proposition}
\label{Pbase} For every~$x \in \Omega$,
\begin{equation} \label{ineq:base}
\begin{aligned}
\mathbb{E} \left[ \max \{ Z_1(x),Z_2(x) \} \right]
& 
= 
\mathbb{E}[ Z_1(x) ] 
\Phi\left(\frac{\diff(x)}{\theta(x)}\right)
+
\mathbb{E}[ Z_2(x) ] 
\Phi\left(\frac{-\diff(x)}{\theta(x)}\right)
+
\theta(x) \phi \left( \frac{\diff(x)}{\theta(x)} \right) \\
& \leq \mathbb{E}\left[Z_1 (x)\right]  
+
u_\theta(x) \frac{1}{\sqrt[]{2\pi}}.
\end{aligned}
\end{equation}
\end{proposition}
\begin{Proof} From the symmetry of the c.d.f. of the normal distribution  follows that
\begin{equation*}
\Phi
\left(\frac{\diff(x)}{\theta(x)}\right)
+ 
\Phi\left(\frac{-\diff(x)}{\theta(x)}\right) = 1.
\end{equation*}
Since $\mathbb{E}\left[Z_1 (x)\right] \geq \mathbb{E}\left[Z_2 (x)\right]$ by assumption, we have
\begin{equation} \label{ineq:base1}
\mathbb{E}\left[Z_1 (x)\right] \geq \Phi\left(\frac{\diff(x)}{\theta(x)}\right)\mathbb{E}\left[Z_1 (x)\right] 
+ 
\Phi\left(\frac{-\diff(x)}{\theta(x)}\right) \mathbb{E}\left[Z_2 (x)\right],
\end{equation}
which constitutes an upper bounding expression for the first two terms of~$\mathbb{E} \left[ \max \{ Z_1(x),Z_2(x) \} \right]$. 

Noting that $\phi\left(\frac{\diff(x)}{\theta(x)}\right) \leq \frac{1}{\sqrt[]{2\pi}}$, and $\theta(x) \leq u_\theta(x)$ by assumption, we have 
\begin{equation}\label{ineq:base2}
u_\theta(x) \frac{1}{\sqrt[]{2\pi}} \geq \theta(x) \phi\left(\frac{\diff(x)}{\theta(x)}\right).  
\end{equation}
Finally, Inequality~\eqref{ineq:base} follows from the addition of Inequality~\eqref{ineq:base1} with Inequality~\eqref{ineq:base2}. $\blacksquare$
\end{Proof}
We propose a discretization approach to obtain $u_\theta(x)$ in Section \ref{discretization1} and present the mathematical model for the baseline RMP in the Appendix. We use this baseline bounding function as a benchmark to measure the improvements in computational performance added by the proposed enhanced bounding techniques in the following sections.  


\subsection{Enhanced Upper-Bounding Function}
\label{sec:Bounds}

Our enhanced approach is based on a joint discretization of $\theta(x)^2$ and $\delta(x)$. The proposed enhanced RMP formulation relies on the following proposition, which is valid for any functions $u_\theta(x), \ l_\theta(x)$, $u_\delta(x)$, and $l_\delta(x)$ such that $0  \leq l_\theta(x) \leq \theta(x) \leq u_\theta(x)$ and $0  \leq l_\delta(x) \leq \delta(x) \leq u_\delta(x)$: 

\begin{proposition}
\label{P2} For every~$x \in \Omega$,
\begin{equation}\label{ineq:ubfunction2}
\mathbb{E}\left[Z_1 (x)\right]  \Phi\left(\frac{u_\delta(x)}{l_\theta(x)}\right)
+ 
\mathbb{E}\left[Z_2 (x)\right] \left(1-\Phi\left(\frac{u_\delta(x)}{l_\theta(x)}\right)\right)
+
u_\theta(x) \phi\left(\frac{l_\delta(x)}{u_\theta(x)}\right)
\geq \mathbb{E} \left[ \max \{ Z_1(x),Z_2(x) \} \right].
\end{equation}
\end{proposition}
\begin{Proof}
This follows because
\begin{equation*}
    \begin{aligned}
        \mathbb{E} \left[ \max \{ Z_1(x),Z_2(x) \} \right]
        & \leq
\mathbb{E}\left[Z_1 (x)\right] \Phi\left(\frac{\delta(x)}{\theta(x)}\right)
+ 
\mathbb{E}\left[Z_2 (x)\right] \left(1-\Phi\left(\frac{\delta(x)}{\theta(x)}\right) \right)
+
u_\theta(x)\phi\left(\frac{\delta(x)}{\theta(x)}\right) \\
& 
\leq 
\mathbb{E}\left[Z_1 (x)\right] \Phi\left(\frac{u_\delta(x)}{l_\theta(x)}\right)
+ 
\mathbb{E}\left[Z_2 (x)\right] \left(1-\Phi\left(\frac{u_\delta(x)}{l_\theta(x)}\right)\right)
+
u_\theta(x) \phi\left(\frac{l_\delta(x)}{u_\theta(x)}\right).
    \end{aligned}
\end{equation*}
The first inequality follows from the symmetry of the c.d.f.  and  because $u_\theta(x) \geq \theta(x)$. 
For the second inequality, first note that, for any constants $a,b$ with $a \geq b$, $a \lambda_1 + b  (1-\lambda_1) \leq a \lambda_2 + b  (1-\lambda_2)$, for every $0 \leq \lambda_1 \leq \lambda_2 \leq 1$.
As the 
c.d.f. is non-decreasing on its domain and $\frac{\delta(x)}{\theta(x)} \leq \frac{u_\delta(x)}{l_\theta(x)}$, and
as 
$\phi\left( y \right)$ is non-increasing for $y \geq 0$  and $\frac{\delta(x)}{\theta(x)} \geq \frac{l_\delta(x)}{u_\theta(x)}$, the result follows.  $\blacksquare$
\end{Proof}

Suppose we are given $d$ intervals $\left\{[\theta^2_q,\theta^2_{q+1}]\right\}_{q=1}^{d}$ and $l$ intervals $\left\{[\delta_h,\delta_{h+1}]\right\}_{h=1}^{l}$, with $\theta(x)^2 \in [\theta^2_1,\theta^2_{d+1}]$ and $\delta(x) \in [ \delta_1 , \delta_{l+1} ]$ for every $x \in \Omega$. 
Furthermore, let $\theta_q$ and $\theta_{q+1}$ denote a lower and upper bound of $\theta(x)$, respectively, for $\theta(x)^2 \in \left[\theta^2_q,\theta^2_{q+1}\right]$.  Using these intervals we construct the following enhanced RMP formulation:
\begin{align}
\max \quad &U+U' \label{RMPObj}\\
\text{s.t. } &u_1 = \sum_{j = 1}^n \mu_j  x_{1,j}; \ u_2 = \sum_{j = 1}^n \mu_j  x_{2,j}; \ u_1 \geq u_2 \label{RMP1}\\
&s = \sum_{i=1}^2\left(\sum_{j=1}^n \sigma^2_j  x_{i,j} + 
2 \sum_{1 \leq j < j' \leq n} \mathrm{cov}(Y_{j} , Y_{j'}) v_{i,j,j'} \right)-2\sum_{j = 1}^{n} \sum_{j' = 1}^n  \mathrm{cov} \left( Y_{j} , Y_{j'}  \right)  r_{j,j'}
\label{RMP2}\\
&v_{i, j,j'} \leq x_{i,j}; \ v_{i,j,j'} \leq x_{i,j'} \hspace{140pt} \forall j,j'\in \set{1,\dots,n}, \ i \in \set{1,2} \label{RMP3}\\
&v_{i, j,j'} \geq x_{i,j}+x_{i,j'}-1 \hspace{154pt} 
\forall j,j' \in \set{1,\dots,n}, \ i \in \set{1,2} \label{RMP4}\\
&r_{j,j'} \leq x_{1,j}; \ r_{j,j'} \leq x_{2,j'} \hspace{149pt} 
\forall j,j' \in \set{1,\dots,n} \label{RMP5}\\
&r_{j,j'} \geq x_{1,j}+x_{2,j'}-1 \hspace{158pt} 
\forall j,j' \in \set{1,\dots,n} \label{RMP6}\\
&  \sum_{q=1}^{d}w_q = 1; \ \sum_{h=1}^{l} y_h = 1; \ s' = \sum_{q=1}^{d} \theta_{q+1} w_q \label{srmp:2} \\
& \theta^2_q w_q \leq s \leq \theta^2_{q+1} + \theta^2_{d+1}(1-w_q) 
\hspace{148pt} q = 1, \ldots, d \label{srmp:3} \\
& \delta_h y_h \leq u_1 - u_2 \leq \delta_{h+1}+ \delta_{l+1}(1-y_h) \hspace{125pt}
h  = 1, \ldots, l  \label{srmp:6} \\
& U \leq 
u_1  \Phi \left( \frac{\delta_{h+1}}{\theta_q} \right)
    +
u_2 \left(1- \Phi \left( \frac{\delta_{h+1}}{\theta_q} \right) \right) 
    +
M 
\left(
    2 - w_q - y_h
\right)
\hspace{10pt} q = 1, \ldots, d, \ h = 1, \ldots, l \label{srmp:7} \\
& U' \leq 
\theta_{q+1} \phi\left(\frac{\delta_{h}}{\theta_{q+1}}\right) 
    +
M 
\left(
    2 - w_q - y_h
\right)
\hspace{108pt} q = 1, \ldots, d, \ h = 1, \ldots, l \label{srmp:7.5} \\
&v \in \set{0,1}^{n \times n \times 2}; \ r \in \set{0,1}^{n \times n}; \ w \in \{0,1\}^{d}; \  y \in \{0,1\}^{l}; \ x \in \Omega. \label{RMP7}
\end{align}
Binary variables~$w_q$, $q = 1, \ldots, d$, and~$y_h$, $h = 1, \ldots, l$, indicate which interval~$\theta(x)^2$ and~$\delta(x)$ belong to, respectively. Variable $u_1$ ($u_2$) denotes $\mathbb{E}\left[Z_1 (x)\right]$ ($\mathbb{E}\left[Z_2 (x)\right]$) and $s$ ($s'$) represents~$\theta(x)^2$ ($u_\theta(x)$). Binary variable $v_{i,j,j'}$  takes a value of 1 iff $x_{i,j}=x_{i,j'}=1$. Similarly, $r_{j,j'}$ equals 1 iff $x_{1,j}=x_{2,j'}=1$. Variable $U$ captures the first two terms of expression \eqref{ineq:ubfunction2} and $U'$ captures the third term. 

The objective function \eqref{RMPObj} maximizes the upper bounding function defined by Proposition \ref{P2}. Constraints \eqref{RMP1} define the $u$-variables according to equation \eqref{expectedVal} and impose the symmetry breaking condition  $u_1 \geq u_2$. Constraint \eqref{RMP2} imposes $s=\theta(x)^2$ as described by equation \eqref{thetaSquared}, where $\sigma^2(Z_1(x))$, $ \sigma^2(Z_2(x))$, and $\mathrm{cov} \left(Z_1(x) , Z_2(x)\right)$ are computed according to equations \eqref{variance} and \eqref{covariance}, respectively. Constraints \eqref{RMP3}--\eqref{RMP6} are the McCormick linearization constraints. Constraints~(\ref{srmp:2}) ensure that exactly one interval is chosen for $\theta(x)^2$ and $\delta(x)$, and set $s'$ equal to the upper bound of $\theta(x)$ for the interval that $\theta(x)^2$ belongs to. Constraints (\ref{srmp:3})--(\ref{srmp:6}) select the right interval for $\theta(x)^2$ and $\delta(x)$. Constraints~(\ref{srmp:7}) are only active for the selected intervals and enforce that~$U$ is bounded by a linear combination of $u_1$ and $u_2$, defined by the evaluation of the c.d.f. at appropriately chosen constants associated with the intervals that~$\theta(x)^2$ and~$\delta(x)$ lie in; where $M$ is a sufficiently large value. Similarly, constraints~(\ref{srmp:7.5}) enforce $U'$ to equal the third term of expression \eqref{ineq:ubfunction2} for the corresponding intervals.    Constraints~\eqref{RMP7} define the domains of the variables appropriately. 


\subsubsection{Discretization of~$\theta(x)^2$\\} \label{discretization1} 
We obtain an upper bound~$\theta^2_{d+1}$ for~$\theta(x)^2$ by solving $\max\{s \ | \ (x,v,r,u_1,u_2,s) \in \Psi\}$, where $\Psi$ is the space defined by constraints \eqref{RMP1}--\eqref{RMP6} and \eqref{RMP7}. This problem is computationally challenging (NP-hard from Theorem~\ref{thm:hard}), so
our strategy consists of solving the resulting MILP for a limited amount of time in order to obtain a relatively refined upper bound.

We define $d$ intervals for $\theta(x)^2$ as follows: $\theta^2_1 = 0, \ \theta^2_2 = 1$,  and $\theta^2_q =\theta^2_{q-1}+ \frac{\theta_{d+1}^2-1}{d-1}$ for $q = 3, \ldots, d+1$. The first interval is different from the others since $\sqrt{a^2} \geq a^2$ for $0\leq a \leq 1$. Note that for any~$x$ in $\Omega$, $0 \leq \theta(x)^2 \leq \theta^2_{d+1}$, 
so $\theta(x)^2$ must belong to $[\theta^2_{q},\theta^2_{q+1}]$ for some 
$q = 1,...,d$.
Given these intervals, upper bounds and lower bounds for~$\theta(x)$ are given by 
\[
\theta_{q+1} = 
\begin{cases}
    1, & q = 1 \\
    \sqrt{\theta^2_{q+1} } , & q = 2, \ldots, d,
\end{cases}
\qquad
\theta_q = 
\begin{cases}
    0, & q = 1 \\
    \sqrt{\theta^2_{q}}, & q = 2, \ldots, d.
\end{cases}
\]

\subsubsection{Discretization of~$\delta(x)$\\} \label{discretization2}
The procedure is analogous to the one described in~\ref{discretization1}. Namely,  we obtain an upper bound~$\delta_{l+1}$ for all $\delta(x)$ by solving problem~$\max\{u_1 - u_2 \ | \ (x,v,r,u_1,u_2,s) \in \Psi \}$ for a limited amount of time. Given~$\delta_{l+1}$, 
we generate $l$ discretization intervals defined by $\delta_h = \frac{h-1}{l} \hat{\delta}$ for $h=1, \ldots l+1$. By construction, $\delta(x)$ belongs to one of the intervals defined by the values~$\delta_h$, as desired.

\subsection{Tightness of the Upper-Bounding Function}\label{sec:ub_tightness}
We investigate now $\Delta(x) = g(x) - \mathbb{E} \left[ \max \{ Z_1(x),Z_2(x) \} \right]$, the difference between the upper bound 
\begin{equation}\label{eq:bound}
g(x) = 
\mathbb{E}\left[Z_1 (x)\right]  \Phi\left(\frac{u_\delta(x)}{l_\theta(x)}\right)
+ 
\mathbb{E}\left[Z_2 (x)\right] \left(1-\Phi\left(\frac{u_\delta(x)}{l_\theta(x)}\right)\right)
+
u_\theta(x) \phi\left(\frac{l_\delta(x)}{u_\theta(x)}\right),
\end{equation}
and the exact expression for~$\mathbb{E} \left[ \max \{ Z_1(x),Z_2(x) \} \right]$, which is given by
\begin{eqnarray*}
\Delta(x) 
&=&
\mathbb{E}\left[Z_1 (x)\right] 
\left( 
\Phi\left(\frac{u_\delta(x)}{l_\theta(x)}\right)
-
\Phi\left(\frac{\delta(x)}{\theta(x)}\right)
\right)
+ 
\mathbb{E}\left[Z_2 (x)\right] \left(
\Phi\left(\frac{\delta(x)}{\theta(x)}\right) - 
\Phi\left(\frac{u_\delta(x)}{l_\theta(x)}\right)
\right) + \\
&& \quad 
u_\theta(x)  \phi\left(\frac{l_\delta(x)}{u_\theta(x)}\right)
-
\theta(x)  \phi\left(\frac{\delta(x)}{\theta(x)}\right) \\
&=&
\delta(x) \left(
\Phi\left(\frac{u_\delta(x)}{l_\theta(x)}\right)
-
\Phi\left(\frac{\delta(x)}{\theta(x)}\right) 
\right) + 
u_\theta(x)  \phi\left(\frac{l_\delta(x)}{u_\theta(x)}\right)
-
\theta(x)  \phi\left(\frac{\delta(x)}{\theta(x)}\right). 
\end{eqnarray*}
The inflection points of~$\Phi$ and~$\phi$ are 0 and $\pm 1$, respectively, so their  first derivatives are bounded by
$\frac{1}{\sqrt{2\pi}}$ and~$\frac{1}{\sqrt{2e\pi}}$, respectively. As
$
\frac{l_\delta(x)}{u_\theta(x)}
\leq
\frac{\delta(x)}{\theta(x)} 
\leq
\frac{u_\delta(x)}{l_\theta(x)}
$,
we have
\[
\Phi\left(
\frac{u_\delta(x)}{l_\theta(x)}
\right)
 \leq
\Phi\left(\frac{\delta(x)}{\theta(x)}
\right)
 + 
\frac{1}{\sqrt{2\pi}}
\left(
\frac{u_\delta(x)}{l_\theta(x)}
-
\frac{\delta(x)}{\theta(x)}
 \right).
\]
Similarly,
\[
\phi\left(
\frac{l_\delta(x)}{u_\theta(x)}
\right)
 \leq
\phi\left(\frac{\delta(x)}{\theta(x)}
\right)
 + 
\frac{1}{\sqrt{2e\pi}}
\left(
\frac{\delta(x)}{\theta(x)}
-
\frac{l_\delta(x)}{u_\theta(x)}
 \right).
\]
Therefore, we have
\begin{eqnarray}\label{deltaBound}
\Delta(x) 
&\leq& 
\frac{\delta(x)}{\sqrt{2\pi}}\left(
\frac{u_\delta(x)}{l_\theta(x)}
-
\frac{\delta(x)}{\theta(x)}
 \right)
 +
 \frac{\theta(x)}{\sqrt{2e\pi}}
\left(
\frac{\delta(x)}{\theta(x)}
-
\frac{l_\delta(x)}{u_\theta(x)}
 \right)
 +
(u_\theta(x) - \theta(x))\phi\left(\frac{\delta(x)}{\theta(x)}\right)  \nonumber \\
&\leq&
\frac{\delta(x)}{\sqrt{2\pi}}\left(
\frac{u_\delta(x)}{l_\theta(x)}
-
\frac{l_\delta(x)}{u_\theta(x)}
 \right)
 +
 \frac{\theta(x)}{\sqrt{2e\pi}}
\left(
\frac{u_\delta(x)}{l_\theta(x)}
-
\frac{l_\delta(x)}{u_\theta(x)}
 \right)
 +
 \frac{u_\theta(x) - \theta(x)}{\sqrt{2e\pi}} \nonumber \\
 &=&
 \left( \frac{\delta(x)}{\sqrt{2\pi}}   + \frac{\theta(x)}{\sqrt{2e\pi}}
 \right)
 \left(
\frac{u_\delta(x)}{l_\theta(x)}
-
\frac{l_\delta(x)}{u_\theta(x)}
 \right)
 +
 \frac{u_\theta(x) - \theta(x)}{\sqrt{2e\pi}}
 \end{eqnarray}
First, note that if~$l_\theta(x) = \theta(x) = u_\theta(x)$ and~$l_\delta(x) = \delta(x) = u_\delta(x)$, $\Delta(x) = 0$, as expected. For the element of the first expression above that does not depend on the discretization, we have for every~$x$
\begin{eqnarray}
\frac{\delta(x)}{\sqrt{2\pi}}   + \frac{\theta(x)}{\sqrt{2e\pi}} 
&\leq&  
\frac{\mathbb{E} \left[ \max \{ Z_1(x),Z_2(x) \} \right]}{\sqrt{2\pi}}
+ 
\frac{\mathbb{E} \left[ \max \{ Z_1(x),Z_2(x) \} \right]}{\sqrt{e}} \nonumber \\
&\leq& 
\frac{\sqrt{2\pi}+\sqrt{e}}{\sqrt{2\pi e}}
\mathbb{E} \left[ \max \{ Z_1(x),Z_2(x) \} \right] \approx 1.005\mathbb{E} \left[ \max \{ Z_1(x),Z_2(x) \} \right]
\end{eqnarray}
Therefore, the quality of~$g(x)$ as an estimator of~$\mathbb{E} \left[ \max \{ Z_1(x),Z_2(x) \} \right]$ depends on~$\left(
\frac{u_\delta(x)}{l_\theta(x)}
-
\frac{l_\delta(x)}{u_\theta(x)} \right)$ and~$u_\theta(x) - \theta(x)$. We leverage this result in order to obtain an approximation guarantee delivered by the enhanced RMP for a special case of our featured machine scheduling application in Section~\ref{sec:scheduling}.

\subsection{Strengthening Inequalities for the RMP}
\label{inequalities}
We show next how estimates on~$\diff(x)$ and~$\theta(x)$ can be coupled via supervalid inequalities (SVIs), which potentially eliminate integer solutions without removing all the optimal ones \citep{IsraeliWood02}. We propose SVIs of the form
\begin{equation}\label{SVIs}
\delta(x) \in \left[\delta_h, \delta_{h+1}\right]  \implies \theta(x) \geq \underline{\theta}^h \quad \forall h=1,\dots,l,
\end{equation}
where 
$\underline{\theta}^h$ establishes a lower bound on~$\theta(x)$ for every~$x$ such that~$\delta(x) \in \left[\delta_h, \delta_{h+1}\right]$. First, we show in Proposition \ref{P3} how to obtain a lower bound~$\underline{\theta}(\delta,z^{LB})$  for~$\theta(x)$ for all solutions~$x$ such that $\delta(x) = \delta$ when~$z^{LB}$ is a known lower bound for the exact problem, which can be obtained using any primal heuristic.
Proposition \ref{P4} extends this result to the case in which $\delta(x)$ belongs to a given interval. Both propositions use value~$\bar{\mu} = \max_{x \in \Omega}\mathbb{E}[Z_1(x)]$, which 
typically can be quickly computed.

\begin{proposition} \label{P3}
Given~$\delta$ and~$z^{LB}$, let~$\Omega\left(\diff,z^{LB}\right)$ be the set of solutions~$x$ such that~$\diff(x)=\diff$ and~$\mathbb{E} \left[ \max \{ Z_1(x),Z_2(x) \} \right] \geq z^{LB}$. A lower bound~$\underline{\theta}(\delta,z^{LB})$ of~$\theta(x)$ for all $x$ in~$\Omega\left(\diff,z^{LB}\right)$ is given~by
\begin{align}
\label{auxProb}
\underline{\theta}(\delta,z^{LB})
= \min_{\theta \geq 0}\left\{ \theta \ | \ \bar{u}+\theta \phi\left(\frac{\delta}{\theta}\right) \geq z^{LB} \right\}.
\end{align}
\end{proposition}

\begin{proposition}
\label{P4}
Given~$\diff_1$, $\diff_2$, and~$z^{LB}$ such that  
$\diff_1 \leq \diff_2$, we have
$\underline{\theta}(\delta_1,z^{LB}) \leq \underline{\theta}(\delta_2,z^{LB})$.
\end{proposition}
\begin{corollary}\label{cor:lb}
If $\delta(x) \in [\delta_h,\delta_{h+1}]$ and 
$\mathbb{E} \left[ \max \{ Z_1(x),Z_2(x) \} \right] \geq z^{LB}$, we have 
$\underline{\theta}\left(\delta_h,z^{LB}\right) \leq \theta(x)$.
\end{corollary}
By setting~$\underline{\theta}^h = \underline{\theta}\left(\delta_h,z^{LB}\right)$, 
we obtain the following valid inequality to RMP:
\[
s \geq \left( \underline{\theta}^h \right)^2 y_h, \quad \quad \forall h=1,\dots,l.
\]
These inequalities tighten the bounds by relating the choices of $h$ for $y_h = 1$ to the best-known solution. Similarly to~$z^{LB}$, each~$\underline{\theta}^h$ is computed at the initialization of Algorithm~\ref{a1}.


\section{Performance Evaluation}
\label{sec:compExperiments}
We present results of an extensive evaluation on synthetic instances to assess the performance of the enhanced algorithm compared to the baseline approach. The models and algorithms are implemented in CPLEX version 12.7.1 \citep{cplex} 
through the Java API.  We utilize Python 3.6 with \texttt{Scikit-learn} functions~\citep{scikit-learn} to generate the random instances. All experiments are conducted on an Intel(R) Xeon(R) CPU E5-1650 v4 at 3.60GHz with 32GB of memory, and we impose a maximum time limit of 10 minutes. Source code and synthetic instances are available upon request.

The problem used for evaluation is the following: \begin{align}
& \min && \mathbb{E} \left[ \max \left\{ Z_1(x),Z_2(x) \right\} \right]  \label{fa:kp} \tag{KP} \\
& \textnormal{s.t.} && Z_i(x) = \sum_{j=1}^n Y_j x_{i,j} && i = 1, 2 \nonumber \\
&&& \sum_{j=1}^n a_j x_{i,j} \leq b_i && i=1,2 \nonumber \\
&&& x_{1,j} + x_{2,j} \leq 1 && j=1, \ldots, n \nonumber \\
&&& x_{i,j} \in \set{0,1} && i=1, 2, \ j = 1, \ldots, n. \nonumber
\end{align}

We generate random problem instances by considering two knapsack constraints defined over a set of $n$ items. Each item has an integer weight $a_j$, which is drawn independently from $\mathcal{U}(1,19)$, and the profit associated to each item, $Y_j$, follows a normal distribution, with mean sampled also from $\mathcal{U}(15,25)$. The right-hand side of the constraints is fixed to~$b_1=b_2=40$. The objective function is to maximize the expected value of the knapsack with maximum profit, where the component random variables $Y_j$ correspond to the uncertain profit of each item. 

To generate variances and correlations for the profits, we use \texttt{Scikit-learn} functions to generate a random, positive semidefinite matrix (PSD), which is then multiplied by a constant factor~$\alpha$. We generate 5 instances for each configuration of~$(n,\alpha)$, where $n \in \{15,20,25,30\}$ and $\alpha$ in~$\{50,100,150,200,250\}$, for a total of 100 instances. The $\alpha$ multiplier allows us to evaluate how sensitive the algorithms are to increasing variance and correlation. 

\begin{figure}[h!]
\centering
\caption{Cumulative distribution plot of performance comparing the baseline and the enhanced optimization algorithms.}
\label{fig:CDPP_Knapsack}
\includegraphics[scale=0.6]{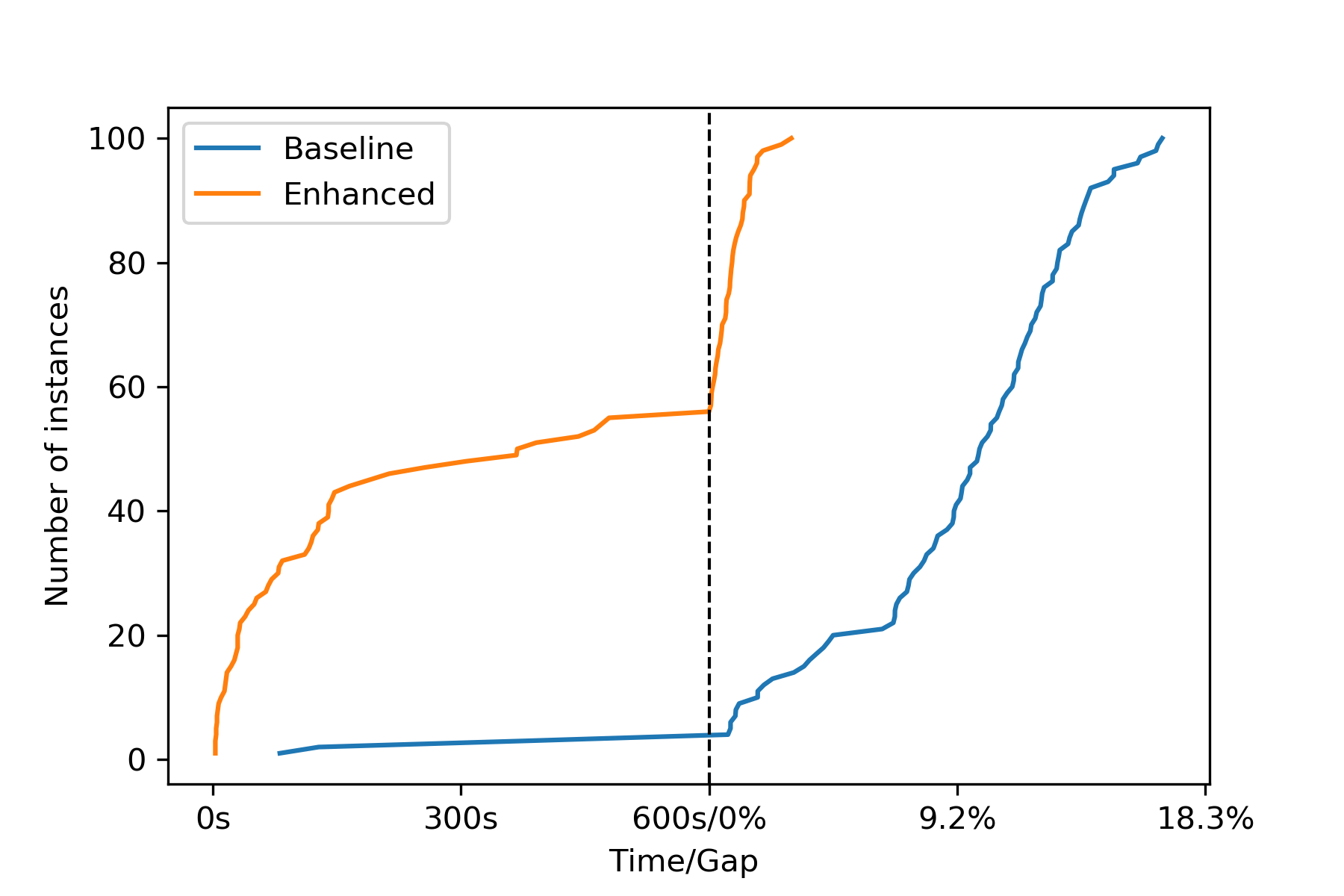}
\end{figure}

Figure~\ref{fig:CDPP_Knapsack} presents a performance plot comparing solution times and optimality gaps for the baseline and enhanced algorithm. On the left half of the graph, the line plots report the number of instances solved by a given time limit, up to 600 seconds.  On the right half of the graph, the line plots report the number of instances solved to within an optimality gap threshold, up to 18.3\%.

Figure~\ref{fig:CDPP_Knapsack} provides strong evidence for (1) the efficacy of the proposed enhancements and (2) the complexity of solving these problems to optimality. The baseline model solves only 3 instances within 10 minutes, and results in significant optimality gaps at time limit for the remaining 96 instances, at 9.3\% on average and up to 18.3\%. With our proposed enhancements, 55 instances are solved within 10 minutes, with 32 solved within 100 seconds, and much smaller optimality gaps for those instances that are not solved to optimality, at 0.9\% on average and 3\% at most. Regarding the number of cuts required to prove optimality, the baseline algorithm generates on average 270.5 cuts, while the enhanced algorithm generates 18.3 cuts, thus demonstrating the superior quality of the upper bounds obtained by the enhanced bounding function.

\renewcommand{\arraystretch}{0.6}
\begin{table}[h!]
\centering
\caption{Tabular comparison of solution times and gaps for baseline and enhanced algorithm.}
\label{tab:baselineVenchanced}
\begin{tabular}{cc|cc|cc}
\multicolumn{2}{c|}{} & \multicolumn{2}{c|}{Baseline} &
\multicolumn{2}{c}{Enhanced} \\
$n$ & $\alpha$ & Time (s)               & Gap (\%) & Time (s)               & Gap (\%) \\ \hline
15	&	50	&	$	254^{(2)}	$	&	1.7	&	$	28^{(5)}	$	&	0.0	\\
15	&	100	&	$	81^{(1)}	$	&	2.6	&	$	44^{(5)}	$	&	0.0	\\
15	&	150	&	$	-^{(0)}	$	&	4.7	&	$	102^{(5)}	$	&	0.0	\\
15	&	200	&	$	-^{(0)}	$	&	4.6	&	$	49^{(5)}	$	&	0.0	\\
15	&	250	&	$	-^{(0)}	$	&	4.8	&	$	28^{(5)}	$	&	0.0	\\ \hline
20	&	50	&	$	-^{(0)}	$	&	7.5	&	$	132^{(5)}	$	&	0.0	\\
20	&	100	&	$	-^{(0)}	$	&	9.4	&	$	87^{(2)}	$	&	0.5	\\
20	&	150	&	$	-^{(0)}	$	&	9.8	&	$	109^{(4)}	$	&	0.2	\\
20	&	200	&	$	-^{(0)}	$	&	10.6	&	$	201^{(4)}	$	&	0.1	\\
20	&	250	&	$	-^{(0)}	$	&	10.8	&	$	111^{(5)}	$	&	0.0	\\ \hline
25	&	50	&	$	-^{(0)}	$	&	10.6	&	$	-^{(0)}	$	&	1.0	\\
25	&	100	&	$	-^{(0)}	$	&	11.3	&	$	-^{(0)}	$	&	0.8	\\
25	&	150	&	$	-^{(0)}	$	&	12.2	&	$	415^{(2)}	$	&	0.6	\\
25	&	200	&	$	-^{(0)}	$	&	12.5	&	$	367^{(1)}	$	&	0.5	\\
25	&	250	&	$	-^{(0)}	$	&	11.9	&	$	213^{(3)}	$	&	0.3	\\ \hline
30	&	50	&	$	-^{(0)}	$	&	10.5	&	$	30^{(1)}	$	&	0.9	\\
30	&	100	&	$	-^{(0)}	$	&	11.0	&	$	71^{(1)}	$	&	0.7	\\
30	&	150	&	$	-^{(0)}	$	&	13.4	&	$	-^{(0)}	$	&	1.4	\\
30	&	200	&	$	-^{(0)}	$	&	13.3	&	$	391^{(1)}	$	&	1.1	\\
30	&	250	&	$	-^{(0)}	$	&	12.7	&	$	306^{(1)}	$	&	0.5	
\end{tabular}
\end{table}
Table~\ref{tab:baselineVenchanced} reports more detailed information concerning the comparison of the baseline and enhanced algorithms.  For each configuration of $(n,\alpha)$ and for both the baseline and enhanced algorithms, we report the average solution time for those instances solved to optimality within 10 minutes, with the number of instances solved to optimality in superscript, and the average percent gap over all instances. This table makes it even more clear that the enhanced algorithm provides significant improvement and suggests that instances become significantly harder to solve for larger values of the variance and correlation.


\section{
Daily Fantasy Sports}
\label{sec:DFF}

In daily fantasy football, contests are arranged based on the starting times of each week's slate of NFL games, and only players in those games are eligible for inclusion on a fantasy roster.  Further, not all players 
are eligible for rosters because the fantasy scoring system only rewards points for specific tasks. Namely, only the quarterbacks (QB), kickers (K), and offensive ``skill position'' players\textemdash wide receivers (WR), running backs (RB), and tight ends (TE)\textemdash are eligible as single players; other players can be
selected collectively as ``team'' defenses (DEF).  

Generally the individual players receive points for gaining yardage and scoring points in the actual game (via touchdowns, extra points, two-point conversions, and field goals) and the team defense earns  points by preventing the opposing team from scoring or by scoring game points itself.

\texttt{DraftKings} is one of the two major DFS providers.  Different types of 
contests are offered on the betting platform, including showdowns, classics, tiers, and others.  We focus exclusively on showdown contests in this application.  Showdown contests only include the players in a single NFL game, and entries consist of six players, regardless of position, with one designated as the captain. The captain costs 1.5 times the normal salary and earns 1.5 times the normal 
points.  Each player may appear no more than once on a given entry, and different contests allow a different number of entries per participant (some as high as 150).  Additionally, each entry requires at least one player from each team in the contest. We focus here on smaller, high-entry-fee contests that limit each participant to two or three entries (we always enter just two). 

DFS contests align directly with the focus of this paper.  First, when assembling a collection of entries, a natural goal for a participant would be to have one of his entries score very high.  This is because of the payout structure\textemdash most of the payout goes to the highest few entries in the competitions, as we will discuss below. This therefore can be modeled as a problem of maximizing the expected value of the higher of the two entries.

\subsection{Problem definition}
We model the two-entry selection problem in a showdown contest as a special case of Problem~\ref{eqn:optProblem}. Let~$n'$ be the number of players, and $n = 2 n'$.  The first~$n'$ players represent standard versions of the players (``flex'' in the lingo used by \texttt{DraftKings}), whereas the next $n'$ represent their ``captain" versions, in a common order. As shown in formulation \eqref{fa:dff}, the uncertain player scores correspond to the random variables in~\ref{eqn:optProblem}. We define the feasible space $\Omega$ via the following constraints.  Let $\mathcal{I}_1$ and $\mathcal{I}_2$ be the sets containing the players from the first and second team, respectively. For $i \in \{1,2\}$ and $j \in [n]$, binary variable~$x_{i,j}$ indicates the selection of player~$j$ for roster~$i$. The following constraints apply to each roster:
\begin{itemize}
    \item Exactly 5 flex players and 1 captain must be selected:
    \[
    \sum\limits_{j=1}^{n'} 
    x_{i,j} = 5, \sum\limits_{j=n'+1}^{n} x_{i,j} = 1, \quad  i = 1, 2.
    \]
    \item The same player cannot be selected both for a flex and the captain positions:
    \[
    x_{i,j} + x_{i,j+n'} \leq 1, \quad i = 1, 2, \ j=1,\dots,n'.
    \]
    \item At least one player from each team appears on each roster:
    \[
    \sum\limits_{j \in \mathcal{I}_i} x_{i,j} \geq 1, \quad i = 1, 2.
    \]   
\end{itemize}

We also limit the set of players under consideration to include only those with an expected score of at least 5 fantasy points.  Players with expectation below this threshold are never selected by our algorithm anyway.  However, as previously mentioned, the assumption of normality in player scores becomes less likely to be rejected as projected scores increase.  

\subsection{Data Sources and Estimation} 
There are several parameters that need to be estimated\textemdash in particular, $\forall j \in [n']$, parameters $\mu_j$ and $\sigma_j$, defining the normal distribution for the points scored by player~$j$; and, $\forall j, j' \in [n']$, the covariance $\rho_{j,j'}$ of the performance of players~$j$ and~$j'$. We used a training set consisting of historical data from four NFL seasons (2014\textendash2017) to estimate those parameters.  We briefly discuss the process for each in turn below.  We then compare our algorithm's performance against a heuristic over 16 competitions in the 2018 season.  More details on the evaluation can be found in the Appendix.

\subsubsection{Expected Value Estimation} Due to the growth of the fantasy sports industry, estimated DFS points for players is the topic of many non-academic articles and websites (e.g., \url{https://rotogrinders.com}).  For the purposes of this paper, we do not generate our own player points projections, but rather use the data from \url{https://fantasydata.com}.  Their projections are consistent and reliable across the years analyzed.

\subsubsection{Variance, Correlation, and Covariance Estimation} We also used the data from
\url{https://fantasydata.com} in order to learn variances, correlations, and covariances for player scores using a nearest-neighbor-like algorithm.  For a player $j$, his variance is estimated as the variance of the actual scores of the 50 players that share a common position with player $j$ in the training set data and whose expected value is as close as possible (in terms of squared difference) to player~$j$.  For example, if a RB is \textit{projected} to score 20.5 fantasy points, we select the 50 RBs in the training set with \textit{projected} fantasy score as close as possible to~20.5 measured by squared difference, and use the \textit{actual} fantasy scores of those 50 players to calculate the variance of that player's fantasy score. 

We estimate the correlation of players $j$ and $j'$ in a similar way to our single-player variance estimates.  Since we are only reporting results for showdown contests, all available player pairs are either on the same or opposing teams.  For teammates, we find the 50 pairs of players and games in the training set for which the players play on the same team, played the same positions as $j$ and $j'$, and have expected values as close as possible to that of $j$ and $j'$, and we use the sample correlation of their actual game scores to estimate their correlation.  For example, if a QB and WR pair are on the same team and are expected to score 30 and 15 points, respectively, we find the 50 instances in the training set of QB and WR teammates with the sum of squared differences from 30 and 15 in expected values as low as possible.  We follow the same process for players on opposing teams.  Finally, a small correction is made to ensure that the covariance matrix is PSD, if needed. 

\subsubsection{Normality Assumption} The assumption for normality of fantasy points production for players is well-grounded.  Using the Shapiro-Wilk test \citep{ShapiroWilk1965}, the null hypothesis of a normal distribution of player scores cannot be rejected for 76\% of the QBs (by far the most valuable position) from the 2016-2018 seasons.  In addition, the null hypothesis of a normal distribution cannot be rejected nearly half the time when considering players from any position over that same time frame with an expected score of at least 10 points.  As the expected score increases, the likelihood of rejecting the assumption of normally distributed actual scores falls across all players.  So while some of the less important positions, which score fewer points, may not follow a normal distribution, we cannot reject the normal distribution for the vast majority of the QBs and in general for the players projected to score a reasonably high number of points.  These and other results regarding the normal distribution of player scores are available from the authors, by request.  

\subsection{Benchmark Heuristic} We used entries generated by a simple heuristic that would select a first entry with maximum expected value and a second entry that differs from the first by at least one player and otherwise again maximizes the expected value.  The set of players available was the same for both our algorithm and the heuristic. Note that this heuristic is similar to the online tools many fantasy participants pay to use, such as \url{https://www.fantasycruncher.com/}.

\subsection{Results} 
Over a collection of 16 contests in the 2018 NFL season, the net realized profit of employing our exact two-entry model using the inputs described above would have been over \$5,000.  On the same collection of contests, employing the heuristic and the same inputs would have resulted in a net loss of over \$4,000. This application therefore highlights how important the joint decision-making over the entries is in order to get a high scoring entry.

We provide detailed results in the Appendix.  One important takeaway is that the two entries selected often score on opposite sides of their expectation, showing how correlation is exploited to elevate the expected score of the maximum.  As an example,
in the game \texttt{Redskins vs. Saints} played on 10-8 the two entries selected had expectation of 104.74 and 102.74, respectively.  Their actual scores were 65.55 and 140.80, where the second entry would have won the competition and resulted in substantial payout.  When we use the heuristic, the two entries selected have expectation 109.51 and 108.88, and scored 78.95 and 114.20.  The score of 114.20 would have resulted in a positive payout, but a very marginal one. This example shows how important optimizing exactly can be in practice.

\section{
Makespan Minimization with Stochastic Processing Times}
\label{sec:scheduling}

We applied our theoretical and algorithmic results to~$P2|p \sim \N(\mu,\Sigma),|\mathbb{E}[C_{max}]$, a machine scheduling problem involving  the minimization of makespan on two parallel machines for jobs with processing times drawn from a multivariate gaussian distribution. We present two theoretical  results for the case where the processing times of the jobs are uncorrelated. First, we show that $P2|p \sim 
\N(\mu,\Sigma),\rho_{j,j'}=0|\mathbb{E}[C_{max}]$
is equivalent to $P2|p_{j} = \mu_j|C_{max}$, the deterministic version of the problem where the means are used as  
processing times.
\begin{theorem}\label{thm:Sched}
$P2|p \sim \N(\mu,\Sigma),\rho_{j,j'}=0|\mathbb{E}[C_{max}]$ is equivalent
to  $P2|p_{i,j} = \mu_{i,j}|C_{max}$.
\end{theorem}
\begin{Proof} Because there is no correlation between the random variables, we have $\theta(x) = 
\sqrt{\sigma^2(Z_1(x))+\sigma^2(Z_2(x))}$. Similarly, we have
$\sigma^2(Z_i) = \sum\limits_{j=1}^n \sigma^2_j  x_{i,j}, i = 1,2$,
and as $x_{1,j} + x_{2,j} = 1$ for each~$j \in [n]$, we have
$\theta(x) = \sqrt{\sum\limits_{j=1}^n \sigma^2_j}$,
i.e., $\theta(x)$ is actually constant and can simply be written as~$\theta$. 
If we set~$\diff_{\theta} = \diff_{\theta}(x) = \frac{
\delta(x)}{\theta}$, we can rewrite
$\mathbb{E} \left[ \max \{ Z_1(x),Z_2(x) \} \right]$ as follows:
\begin{align*}\label{eq:objective}
& \mathbb{E} \left[ \max \{ Z_1(x),Z_2(x) \} \right]
= 
\mathbb{E} 
\left[ 
	Z_1 (x)
\right]
\Phi
\left(
\diff_{\theta}
\right)
+ 
\mathbb{E} 
\left[ 
	Z_2 (x)
\right]
\Phi
\left(
    	\diff_{\theta}
\right) +
\theta
\phi
\left(
    	\diff_{\theta}
\right).
\end{align*}
The c.d.f. of the standard normal distribution can be written as
\[
\Phi(x) = 
\frac{1}{2} + \frac{1}{\sqrt{2\pi}}e^{\frac{-x^2}{2}}\left[ x + \frac{x^3}{3} + \frac{x^5}{5\cdot 3} + \ldots + \frac{x^{2n+1}}{(2n+1)!!} + \ldots \right],
\]
where $n!! = n (n-2) (n-4)\ldots(((n-1) \mod 2)+1)$ is the double factorial of $n$. Therefore,
\begin{eqnarray*}
\mathbb{E} \left[ \max \{ Z_1(x),Z_2(x) \} \right] &=& 
\mathbb{E} \left[ Z_1(x) \right] \left( \frac{1}{2} + \frac{1}{\sqrt{2\pi}}e^{\frac{-\diff_{\theta}^2}{2}}\left[\diff_{\theta} + \frac{\diff_{\theta}^3}{3!!} + \frac{\diff_{\theta}^5}{5!!} + \ldots  \right]\right) + \\
&& 
\mathbb{E} \left[ Z_2(x) \right]\left( \frac{1}{2} - \frac{1}{\sqrt{2\pi}}e^{\frac{-\diff_{\theta}^2}{2}}\left[\diff_{\theta} + \frac{\diff_{\theta}^3}{3!!} + \frac{\diff_{\theta}^5}{5!!} + \ldots  \right]\right) + 
\frac{\theta}{\sqrt{2\pi}}e^{\frac{-\diff_{\theta}^2}{2}} \\
&=&
\frac{\mathbb{E} \left[ Z_1(x) \right]+\mathbb{E} \left[ Z_2(x) \right]}{2} + \left( \frac{
\theta \diff_{\theta}}{\sqrt{2\pi}}e^{\frac{-\diff_{\theta}^2}{2}}\left[\diff_{\theta} + \frac{\diff_{\theta}^3}{3!!} + \frac{\diff_{\theta}^5}{5!!} + \ldots  \right]\right) + \frac{\theta}{\sqrt{2\pi}}e^{\frac{-\diff_{\theta}^2}{2}} \\
&=&
\frac{\mathbb{E} \left[ Z_1(x) \right]+\mathbb{E} \left[ Z_2(x) \right]}{2} + 
\left( \frac{\theta}{\sqrt{2\pi}}e^{\frac{-\diff_{\theta}^2}{2}}\left[1 + \diff_{\theta}^2 + \frac{\diff_{\theta}^4}{3!!} + \frac{\diff_{\theta}^6}{5!!} + \ldots  \right]\right). 
\end{eqnarray*}
By taking the first derivative in~$\diff_{\theta}$, we obtain
\begin{eqnarray*}
\frac{d \mathbb{E} \left[ \max \{ Z_1(x),Z_2(x) \} \right] }{d\diff_{\theta}} &=&
\left( \frac{\theta}{\sqrt{2\pi}}e^{\frac{-\diff_{\theta}^2}{2}}\left[2\diff_{\theta} + \frac{4\diff_{\theta}^3}{3!!} + \frac{6\diff_{\theta}^5}{5!!} + \ldots  \right]\right)
-
\left( \frac{\theta}{\sqrt{2\pi}}e^{\frac{-\diff_{\theta}^2}{2}}\left[\diff_{\theta}  + \diff_{\theta}^3 + \frac{\diff_{\theta}^5}{3!!} + \frac{\diff_{\theta}^7}{5!!} + \ldots  \right]\right) \\
&=&
\frac{\theta}{\sqrt{2\pi}}e^{\frac{-\diff_{\theta}^2}{2}}\left[\diff_{\theta} + \frac{\diff_{\theta}^3}{3!!} + \frac{\diff_{\theta}^5}{5!!} + \ldots  \right],
\end{eqnarray*}
which is strictly positive for~$\diff_{\theta} \geq 0$. Therefore, we conclude that~$\mathbb{E} \left[ \max \{ Z_1(x),Z_2(x) \} \right]$ attains its minimum (maximum) for every~$x$ such that~$\diff_{\theta}(x)$ is minimum (maximum).  $\blacksquare$
\end{Proof}
This framing of $P2|p \sim \N(\mu,\Sigma),\rho_{i,j}=0|\mathbb{E}[C_{max}]$ allows one to see that the minimization of makespan is equivalent to maximization of the idlest machine's load. Similar arguments hold if one wishes to optimize~$\mathbb{E} \left[ \min \{ Z_1(x),Z_2(x) \} \right]$; in the machine scheduling setting, this problem consists of minimizing the load of the idlest machine (or maximizing the load of the busiest machine), which admits a trivial optimal solution, where all jobs are assigned to one machine. 

The next result leverages our results in Section~\ref{sec:ub_tightness} involving the tightness of the upper-bound function to shows that RMP  delivers solutions with constant-factor approximation guarantees for~$P2|p \sim \N(\mu,\Sigma),\rho_{i,j}=0|\mathbb{E}[C_{max}]$.

\begin{theorem}
RMP delivers a 2.005-approximation for~$P2|p \sim \N(\mu,\Sigma),\rho_{j,j'}=0|\mathbb{E}[C_{max}]$.
\end{theorem} 
\begin{Proof} As $\theta(x) = \theta$ if variables are uncorrelated, the second term in~\eqref{deltaBound} vanishes and we have
 \[
 \Delta(x)
 =
 \left( \frac{\delta(x)}{\sqrt{2\pi}}   + \frac{\theta}{\sqrt{2e\pi}}
 \right)
 \left(
\frac{u_\delta(x) - l_\delta(x)}{\theta}
 \right).
 \]
 Also, because variances play no role in this setting (from Theorem~\ref{thm:Sched}), one may set all discretization intervals of~$\delta$ to have the same length and scale the variances such that $\theta = u_\delta(x) - l_\delta(x)$ for every~$x$. Therefore,
$ g(x) \leq 
\mathbb{E} \left[ \max \{ Z_1(x),Z_2(x) \} \right]\left(1 + 
\frac{\sqrt{2\pi}+\sqrt{e}}{\sqrt{2\pi e}}\right)
 \approx 2.005\mathbb{E} \left[ \max \{ Z_1(x),Z_2(x) \} \right]$.
$\blacksquare$
\end{Proof}
 
We also conducted a computational evaluation of our algorithm using synthetic instances of~$P2|p \sim \N(\mu,\Sigma)|\mathbb{E}[C_{max}]$ that were generated using the same procedures as~\cite{ranjbar2012two} and~\cite{stec2019scheduling}; details about the generation of these instances are presented in the Appendix (see Section~\ref{sec:experiments_makespan}).  Overall, our algorithm had a solid performance, achieving an average optimality gap of approximately 0.12\% within the time limit of 10 minutes; in particular, 98 out of 180 instances were solved to optimality. In a comparison with a deterministic heuristic that ignores correlations and uses the averages as the execution times of the jobs, the average improvements are roughly~1.7\%; slightly larger differences were observed when the processing times had higher variances.

\section{Conclusion}
\label{sec:conclusion}
We investigate a class of challenging discrete stochastic optimization problems for which the objective function is given as the expected value of the maximum of two functions of the component random variables of a multivariate Gaussian distribution. We show that our problem is NP-hard and provide two real-world applications that can be modeled within our settings. 

From a computational perspective, the main difficulty for solving these problems comes from the highly nonlinear expression describing the objective function, which contains the evaluation of both the c.d.f and p.d.f of a standard normal distribution with arguments given as functions of the decision variables. We propose an exact cutting-plane algorithm based on a linear function that provides upper bounds on the nonlinear objective. We investigate strengthening techniques for the bounding function and our computational results show a considerable improvement in performance as a result of the proposed techniques. 

For the featured applications, computational results show that our algorithm provides a clear advantage over deterministic heuristics that do not consider correlations and covariances. We also present results of theoretical relevance for an stochastic makespan minimization problem with two machines. We prove the equivalence of $2|p_{i,j} \sim \N(\mu_{i,j},\sigma_{i,j}^2),\rho_{j,j'}=0|\mathbb{E}[C_{max}]$ and its deterministic counterpart and show that optimizing over our bounding function gives a 2.005-approximation for the stochastic version. Real-world settings of the problem, such as DFS, allow for scenarios where three or more functions (entries, in this case) can be selected. The resulting problems are challenging from both a mathematical and computational perspective, and investigating them is an exciting possibility for future work.

\bibliographystyle{plainnat} 
\bibliography{dffbib.bib} 

\clearpage 


\appendix

\section{Daily Fantasy Sports Estimation Details and Details of Results}
\label{sec:DFFParamEst}

Further details of Section~\ref{sec:DFF} are shown below.

\subsection{Expected Value Estimation}
\label{sec:EVEst}
The community of DFS participants puts tremendous resources into calculating reliable and accurate estimates for how a player will perform.  Although there is potentially room to improve upon published estimates, for this paper we use estimated player data from \url{https://fantasydata.com}, which is updated frequently for upcoming games and contains historical NFL game data since 2014, including projected fantasy points for players, actual fantasy points for players, and player salaries on the \texttt{DraftKings} platform, all for a 
monthly fee.  

The projected points estimates used for this paper are the final \url{https://fantasydata.com} estimates set just before game time (changes leading up to game time 
may occur due to weather, injury updates, or a myriad of other reasons).   This data therefore provides $\mu_j$, for $j=1, \ldots p'$, and, because $\mathbb{E} [1.5  X] = 1.5  \mathbb{E}[X]$, we have $\mu_j = 1.5 \mu_{j-p'}$, for $j \geq p'+1$.  These \url{https://fantasydata.com} estimates are consistent, and therefore reliable for future use, from the 2014-2017 seasons.  We fail to reject the null hypothesis that the difference between the \url{https://fantasydata.com} estimates and the actual fantasy scores across all players (with projected points above 5) and all games from 2014-2017 is equal to 0 (the 95\% confidence interval on those 11,785 observations is [-0.211, 0.071]).  



\subsection{Variance, Correlation, and Covariance Estimation}
\label{sec:CorEst}

Due to the nature of the sport itself and the way fantasy scoring works, players on the same (or opposing) teams often have correlated scoring.  For example, if a QB throws a touchdown to a WR, both players receive a number of fantasy points (and the opposing DEF would possibly lose some).  But since that QB had to have completed the pass to someone, fantasy points for individual players rarely exist in a vacuum.  Rather, we would expect that players at certain positions would have significantly correlated scores with teammates, and even opponents, at other positions.  Therefore, given the heavily skewed payoff structure of most \texttt{DraftKings} showdown contests, when trying to maximize the expected value of the maximum score of the entries, participants must take into account these correlations.  The DFS betting community is well aware of this strategy of choosing teammates (or opponents) whose scores should correlate by design; the strategy is termed ``stacking."  It should also be noted that opponents may see their actual scores correlate \textit{even if the players are not on the field at the same time} due to the nature of the game of football.  For example, the actual fantasy scores of the QBs on opposing teams (who are \textit{never} on the field at the same time) are positively correlated because higher-scoring NFL games tend to generate more fantasy points for the QBs.  As one QB scores fantasy points (as his on-field team scores lots of points in the game itself), the other team tends to throw more to try and catch up, and this tends to raise the opponent QB's fantasy scores as a result.

Once the correlations for players on the same and opposing teams have been estimated as previously described, one can then set the correlation to 0 if the significance of the Pearson's correlation test (\cite{benesty2009pearson}) is above some threshold~$\eta$.  Further, there are instances in which the correlation identifies values that are incompatible because of the way the estimation is done and the actual scores achieved by certain player pairs.  Formally, the estimated covariance matrix~$\Sigma$ for a multivariate normal distribution must be PSD.   Correlation estimates, multiplied by the corresponding standard deviations, are found only to provide an estimate for $\Sigma$, and the estimation procedure previously described often yields correlations for which $\Sigma$ is not PSD.  There are several ways to correct for this.  One is to use packages available in common statistical software to adjust the covariances, like the function \texttt{cov\_nearest} in Python's \texttt{statsmodels} module~(\cite{seabold2010statsmodels}). Another is to dampen the correlations by a fixed constant (i.e., set $\rho_{j,j'} := \frac{\rho_{j,j'}}{\chi}$, with $\chi > 0$, for all~$j$,$j'$). Yet another way is to apply the PSD condition and, whenever a solution is found for which $\theta < 0$, redefine $\theta := 0$.  We tried several procedures, though our reported results use the \texttt{cov\_nearest} function, using $\eta = 0.25$ and \texttt{cov\_nearest} for correcting $\Sigma$; i.e., setting to 0 any correlation whose $p$-value is above~0.25.

\subsection{Results} 

In Table~\ref{tab:betting cov2}, we report results from actual competitions on the \texttt{DraftKings} platform that would have been obtained on 16 showdown contests from the 2018 season by our exact algorithm with the parameters estimated as described in the Appendix.  In Table~\ref{tab:betting heur}, we report the same information for the heuristic.
These 16 specific contests were selected because we have actual output from real contests for these competitions. They include 2-entry and 3-entry competitions.

The content of the first seven columns of Tables~\ref{tab:betting cov2} and~\ref{tab:betting heur} is 
as follows.
The first column indicates the game and the date when it took place.  The second indicates the total number of entries in the contest.  The third indicates the price paid by participants per entry.  The fourth and fifth indicate the results obtained by showing the winnings and the profit that the entries identified by the algorithm would have yielded, respectively.  Columns six and seven indicate the best score and the minimum score for entries in the money (i.e., those obtaining any payout from the competition).  Columns eight and nine contain the expected values of the first and the second entry obtained, respectively.  Finally, columns ten and eleven represent the actual scores of the selected entries.  The values in columns ten and eleven are \textbf{bolded} if the entry would have finished in a position receiving payout
in that particular contest.
The total entry fees, winnings, and profit for each method are computed in the last row of each table.  Note, however, that the money won shown in Tables~\ref{tab:betting cov2} and~\ref{tab:betting heur} is an approximation in the sense that if our two entries from the algorithms in this paper were \textit{actually} included in the contests, two other entries that were in the contest necessarily would have been excluded, since these contests filled up to capacity.   If, for example, the two ``eliminated" entries were top-scoring ones, the payoffs for the remaining entries, and ours, would necessarily be the same or larger.  There is no way to account for this game of elimination in a reasonable way, so we report results as if our entries were simply added to the contest without changing any other details about other entries and the payout structure.  

There are two initial takeaways 
from Tables~\ref{tab:betting cov2} and~\ref{tab:betting heur}.  First, the contests differ in cost and number of entries considerably, though the most common version of this contest features 100 total entries at a cost of \$444 per entry.  Second, the scores needed to win, or even place in the money, differ tremendously by game. Namely, higher-scoring NFL games result in more fantasy points for the involved players; for instance, in the game
\texttt{Vikings vs. Rams} played on 9-27, which finished with a score of 38-31, the \textit{minimum payout score} would have been the \textit{winner score} in every single other contest shown.  Because of the vast discrepancy in contest results, it is difficult to say \textit{a priori} what score is going to be necessary to win or place in a given contest.  The comparison of overall results across the methods reveals the power of our algorithm, so the more important takeaway is the following:   
our 
exact algorithm would result in a positive return of over 50\% and the heuristic in a loss of nearly 50\%.

In nearly all contests, across both specifications, the winner score is significantly higher than the expected value scores of our entries 
(the only major exception being the contest featuring the game \texttt{Falcons vs. Saints} played on 11-22).  The players on the winning entry, as a whole, significantly outperform their projections in nearly every contest.  Something akin to the heuristic method is what we expect many participants in these contests are using as their DFS strategy.  Although it unsurprisingly outperforms our entries in \textit{expected value}, its actual results are inferior.
As we have observed in our experiments with 
the other featured applications, in general, the maximum expected value entry does not necessarily belong to an optimal configuration of entries, and it also will not necessarily result in a payout, so bettors need to exploit joint decision-making across multiple entries in order to elevate their actual scores above the threshold needed to win money.

Table~\ref{tab:monster evals 2} further illustrates the power of our algorithm for~\ref{eqn:optProblem}.  We again compare to the heuristic by evaluating Expression~\eqref{eq:objective} 
using
the solution that the heuristic obtains with the same covariance matrix we use.  In this case, ``EV" represents the single highest expected value entry, ``AV" represents the single highest actual score of the entries, and ``OF", 
the evaluation of the objective function at the solution obtained.  Notice that while EV for the heuristic is in most cases higher than EV for the exact approach, the evaluation of the objective function is considerably higher for the exact method in all the instances.
The last row of the table considers the average values above it. Over these contests, the average EV is 1.19 points higher for the heuristic but the average objective function is 6.68 points higher for the exact method, once again showing that our approach can exploit correlation in a way that a common heuristic cannot. 

\newpage

\begin{landscape}
\begin{table}[]
\caption{Results from real-world betting scenarios using ~\ref{eqn:optProblem} and corrected covariance}
\label{tab:betting cov2}
\footnotesize
\begin{tabular}{|c|c|c|c|c|c|c|c|c|c|c|}
\hline
Game & Entrants & Entry Fee & Winnings & Profit & Winner Score & Min. Pay Score & EV E1 & EV E2 & AV E1 & AV E2 \\
\hline
Jets vs. Browns on 9-20 & 100 & \$666 & \$0 & -\$666 & 96.51 & 67.96 & 77.97 & 76.55 & 64.07 & 42.02 \\
Seahawks vs. Cowboys on 9-23 & 1,189 & \$8 & \$0 & -\$8 & 109.78 & 76.78 & 79.63 & 77.97 & 75.22 & 61.57 \\
Vikings vs. Rams on 9-27 & 83 & \$800 & \$0 & -\$800 & 205.28 & 175.13 & 98.50 & 95.90 & 159.28 & 150.13 \\ 
Ravens vs. Steelers on 9-30 & 294 & \$40 & \$0 & -\$40 & 111.01 & 96.73 & 93.40 & 89.16 & 84.66 & 86.17 \\ 
Broncos vs. Chiefs on 10-1 & 69 & \$800 & \$0 & -\$800 & 110.36 & 90.69 & 106.88 & 106.46 & 79.69 & 71.50 \\ 
Colts vs. Patriots on 10-4 & 100 & \$888 & \$700 & -\$188 & 148.09 & 131.85 & 93.13 & 92.75 & 114.19 & \textbf{132.34} \\ 
Redskins vs. Saints on 10-8 & 100 & \$888 & \$10,000 & \$9,112 & 134.97 & 99.67 & 104.74 & 102.74 & 65.55 & \textbf{140.80} \\ 
Giants vs. Falcons on 10-22 & 100 & \$888 & \$2,700 & \$1,812 & 143.56 & 120.70 & 103.03 & 101.08 & \textbf{133.76} & \textbf{123.26} \\ 
Raiders vs.49ers on 11-1 & 100 & \$888 & \$1,200 & \$312 & 95.81 & 78.61 & 77.54 & 77.35 & \textbf{86.72} & 43.61 \\ 
Panthers vs. Steelers on 11-8 & 100 & \$888 & \$0 & -\$888 & 149.93 & 118.70 & 103.03 & 99.76 & 107.00 & 115.49 \\ 
Giants vs. 49ers on 11-12 & 100 & \$888 & \$0 & -\$888 & 120.60 & 95.18 & 91.00 & 90.61 & 81.52 & 81.37 \\ 
Packers vs. Seahawks on 11-15 & 100 & \$888 & \$0 & -\$888 & 126.53 & 114.43 & 91.12 & 88.11 & 107.43 & 96.68 \\ 
Steelers vs. Jaguars on 11-18 & 70 & \$66 & \$0 & -\$66 & 115.86 & 96.49 & 92.96 & 88.87 & 94.45 & 81.31 \\
Falcons vs. Saints on 11-22 & 151 & \$150 & \$450 & \$300 & 122.22 & 106.52 & 130.90 & 129.62 & \textbf{110.94} & \textbf{114.12} \\
Redskins vs. Eagles on 12-3 & 100 & \$888 & \$0 & -\$888 & 111.36 & 84.29 & 88.60 & 83.49 & 82.74 & 64.14 \\
Jaguars vs. Titans on 12-6 & 882 & \$40 & \$0 & -\$40 & 134.60 & 74.83 & 78.30 & 77.25 & 54.93 & 47.83 \\ 
\hline
TOTALS & & \$9,674 & \$15,050 & \$5,376 & & & & & & \\
\hline
\end{tabular}
\end{table}
\end{landscape}

\newpage

\begin{landscape}
\begin{table}[]
\caption{Results from real-world betting scenarios using heuristic}
\label{tab:betting heur}
\footnotesize
\begin{tabular}{|c|c|c|c|c|c|c|c|c|c|c|}
\hline
Game & Entrants & Entry Fee & Winnings & Profit & Winner Score & Min. Pay Score & EV E1 & EV E2 & AV E1 & AV E2 \\
\hline
Jets vs. Browns on 9-20 & 100 & \$666 & \$0 & -\$666 & 96.51 & 67.96 & 78.79 & 78.50 & 63.82 & 60.81 \\
Seahawks vs. Cowboys on 9-23 & 1,189 & \$8 & \$0 & -\$8 & 109.78 & 76.78 & 80.93 & 80.91 & 64.20 & 66.70 \\
Vikings vs. Rams on 9-27 & 83 & \$800 & \$0 & -\$800 & 205.28 & 175.13 & 98.50 & 97.91 & 159.28 & 147.68 \\ 
Ravens vs. Steelers on 9-30 & 294 & \$40 & \$0 & -\$40 & 111.01 & 96.73 & 94.07 & 93.40 & 73.80 & 84.66 \\ 
Broncos vs. Chiefs on 10-1 & 69 & \$800 & \$1,200 & \$400 & 110.36 & 90.69 & 109.91 & 109.85 & \textbf{99.60} & 72.75 \\ 
Colts vs. Patriots on 10-4 & 100 & \$888 & \$0 & -\$888 & 148.09 & 131.85 & 94.25 & 93.59 & 119.44 & 112.11 \\ 
Redskins vs. Saints on 10-8 & 100 & \$888 & \$1,000 & \$112 & 134.97 & 99.67 & 109.51 & 108.88 & 78.95 & \textbf{114.20} \\ 
Giants vs. Falcons on 10-22 & 100 & \$888 & \$2,700 & \$1,812 & 143.56 & 120.70 & 103.03 & 103.03 & \textbf{133.76} & \textbf{123.96} \\ 
Raiders vs.49ers on 11-1 & 100 & \$888 & \$0 & -\$888 & 95.81 & 78.61 & 80.21 & 80.21 & 71.49 & 76.79 \\ 
Panthers vs. Steelers on 11-8 & 100 & \$888 & \$0 & -\$888 & 149.93 & 118.70 & 103.93 & 103.59 & 118.50 & 99.00 \\ 
Giants vs. 49ers on 11-12 & 100 & \$888 & \$0 & -\$888 & 120.60 & 95.18 & 93.45 & 93.30 & 80.72 & 88.72 \\ 
Packers vs. Seahawks on 11-15 & 100 & \$888 & \$0 & -\$888 & 126.53 & 114.43 & 92.54 & 92.26 & 105.43 & 105.08 \\ 
Steelers vs. Jaguars on 11-18 & 70 & \$66 & \$0 & -\$66 & 115.86 & 96.49 & 95.26 & 94.92 & 80.25 & 82.67 \\
Falcons vs. Saints on 11-22 & 151 & \$150 & \$400 & \$250 & 122.22 & 106.52 & 130.90 & 130.02 & \textbf{110.94} & \textbf{113.06} \\
Redskins vs. Eagles on 12-3 & 100 & \$888 & \$0 & -\$888 & 111.36 & 84.29 & 88.70 & 88.60 & 78.94 & 82.74 \\
Jaguars vs. Titans on 12-6 & 882 & \$40 & \$0 & -\$40 & 134.60 & 74.83 & 79.53 & 79.38 & 49.17 & 36.73 \\ 
\hline
TOTALS & & \$9,674 & \$5,300 & -\$4,374 & & & & & & \\
\hline
\end{tabular}
\end{table}
\end{landscape}

\newpage

\begin{table}[]
\centering
\caption{Comparing the exact approach  with the benchmark heuristic on 16 DFS competitions.}
\label{tab:monster evals 2}
\footnotesize
\begin{tabular}{|c|c|c|c|c|c|c|}
\hline
& \multicolumn{3}{c|}{Exact} & \multicolumn{3}{c|}{Heuristic} \\
\hline
Game & EV & OF & AV & EV & OF & AV \\
\hline
Jets vs. Browns on 9-20 & 77.97 & 93.30 & 64.07 & 78.79 & 88.35 & 63.82 \\
Seahawks vs. Cowboys on 9-23 & 79.63 & 96.81 & 75.22 & 80.93 & 88.69 & 66.70 \\
Vikings vs. Rams on 9-27 & 98.50 & 115.53 & 159.28 & 98.50 & 105.93 & 159.28 \\ 
Ravens vs. Steelers on 9-30 & 93.40 & 110.12 & 86.17 & 94.07 & 109.82 & 84.66 \\ 
Broncos vs. Chiefs on 10-1 & 106.88 & 133.03 & 79.69 & 109.91 & 125.27 & 99.60 \\ 
Colts vs. Patriots on 10-4 & 93.13 & 107.85 & 132.34 & 94.25 & 107.29 & 119.44 \\ 
Redskins vs. Saints on 10-8 & 104.74 & 131.43 & 140.80 & 109.51 & 119.91 & 114.20 \\ 
Giants vs. Falcons on 10-22 & 103.03 & 119.65 & 133.76 & 103.03 & 107.43 & 133.76 \\ 
Raiders vs.49ers on 11-1 & 77.54 & 94.74 & 86.72 & 80.21 & 85.07 & 76.79 \\ 
Panthers vs. Steelers on 11-8 & 103.03 & 118.79 & 115.49 & 103.93 & 108.71 & 118.50 \\ 
Giants vs. 49ers on 11-12 & 91.00 & 109.97 & 81.52 & 93.45 & 96.02 & 88.72 \\ 
Packers vs. Seahawks on 11-15 & 91.12 & 110.23 & 107.43 & 92.54 & 102.43 & 105.43 \\ 
Steelers vs. Jaguars on 11-18 & 92.96 & 110.96 & 94.45 & 95.26 & 109.91 & 82.67 \\
Falcons vs. Saints on 11-22 & 130.90 & 146.04 & 114.12 & 130.90 & 145.32 & 113.06 \\
Redskins vs. Eagles on 12-3 & 88.60 & 102.29 & 82.74 & 88.70 & 94.51 & 82.74 \\
Jaguars vs. Titans on 12-6 & 78.30 & 93.68 & 54.93 & 79.53 & 92.93 & 49.17 \\ 
\hline
AVERAGES & 94.65 & 112.15 & 100.09 & 95.84 & 105.47 & 97.41 \\
\hline
\end{tabular}
\end{table}

\section{Instances used in the Experiments of the Makespan Minimization Application}\label{sec:experiments_makespan}

We evaluate our algorithms for $P2|p \sim \N(\mu,\Sigma)|\mathbb{E}[C_{max}]$ 
using a set of synthetic instances that follow the same generation procedure  proposed by~\cite{ranjbar2012two} and~\cite{stec2019scheduling}. These papers consider a related objective of maximizing the number of jobs completed within a pre-specified time limit and assume pairwise independence between the processing time of the jobs. 

We generated instances with $n = 15, 20, 25$, and $30$, where~$n$ denotes the number of jobs. The processing time of each job~$j$ follows a normal distribution~$\N(\mu_j,\sigma_j^2)$, where~$\mu_j$ is a positive value drawn from~$\N(20,9)$ and $\sigma^2_j$ is drawn from $\mathcal{U}(0,0.1 \cdot \mu_j^2 \cdot \eta)$ for~$\eta \in \{0.25,0.5,0.75\}$.  To generate covariance matrices, for each job, we assign uniformly at random a number from~$\{1,2,3\}$,  representing one of three clusters, and we assume that the execution time of jobs assigned to the same cluster have a correlation of one; otherwise, jobs in different clusters have correlation of zero. These correlation matrices are suitable for scenarios where processing times are strongly affected by features shared by several jobs; for example, processing times of computational experiments originating from a same batch are expected to be highly correlated. We generated 5 instances for each combination of values of~$n$ and~$\eta$, so our data set for this application consists of 180 instances. 

\section{Adapting the Optimization Algorithm for Minimization}
\label{app:opti}
In the minimization case, the objective function is still given by expression \eqref{eq:objective}; however, in contrast to the maximization case, our RMP should yield lower bounds on the objective. Proposition \ref{P1-2} presents our proposed enhanced lower-bounding function $g$.
\begin{proposition}
\label{P1-2} For every~$x \in \Omega$,
\begin{equation}
\mathbb{E}\left[Z_1 (x)\right]  \Phi\left(\frac{l_\delta(x)}{u_\theta(x)}\right)
+ 
\mathbb{E}\left[Z_2 (x)\right] \left(1-\Phi\left(\frac{l_\delta(x)}{u_\theta(x)}\right)\right)
+
l_\theta(x) \phi\left(\frac{u_\delta(x)}{l_\theta(x)}\right)
\leq \mathbb{E} \left[ \max \{ Z_1(x),Z_2(x) \} \right].
\end{equation}
\end{proposition}
\begin{Proof}
The proof follows the same arguments as the proof of Proposition \ref{P2} . $\blacksquare$
\end{Proof}

Using this lower bounding function it is straightforward to modify the RMP to yield lower bounds. The modified cutting algorithm for the minimization case is as follows: 
\begin{algorithm}[H]
\begin{algorithmic}[1]

\State Set \textit{LB}$ = -\infty$, \textit{UB}$ = \infty$, $\mathcal{C} = \emptyset$, and incumbent solution $\bar{x} = 0$.

\State Optimize RMP$(\mathcal{C})$ to obtain~$\hat{x}$;  if the problem is infeasible, go to Step 6. 

\State 
Set \textit{LB} = $g(\hat{x})$

\State 
If $\mathbb{E} \left[ \max \{ Z_1(x),Z_2(x) \} \right]$ $<$ \textit{UB}, set \textit{UB} = $\mathbb{E} \left[ \max \{ Z_1(x),Z_2(x) \} \right]$, and update incumbent $\bar{x} = \hat{x}$. 

\State If \textit{LB} = \textit{UB}, go to Step 6.  Otherwise, set~$\mathcal{C} = \mathcal{C} \cup \{c(\hat{x})\}$ and 
return to Step 2.

\State If \textit{UB} = $\infty$, 
original problem is infeasible. Otherwise,
terminate with optimal solution 
$\bar{x}$.
\end{algorithmic}
\caption{A Cutting-Plane Algorithm for Minimization Problems}
\label{a2}
\end{algorithm}

\section{Formulation for the Baseline Approach}
Suppose we are given $d$ intervals $\left\{[\theta^2_q,\theta^2_{q+1}]\right\}_{q=1}^{d}$ with $\theta(x)^2 \in [\theta^2_1,\theta^2_{d+1}]$ for every $x \in \Omega$. 
Furthermore, let $\theta_q$ and $\theta_{q+1}$ denote a lower and upper bound of $\theta(x)$, respectively, for $\theta(x)^2 \in [\theta^2_q,\theta^2_{q+1}]$.  Using these intervals we construct the following RMP formulation, where binary variables~$w_q$, $q = 1, \ldots, d$ indicate which interval~$\theta(x)^2$. Variable $u_1$ ($u_2$) denotes $\mathbb{E}\left[Z_1 (x)\right]$ ($\mathbb{E}\left[Z_2 (x)\right]$) and $s$ ($s'$) represents~$\theta(x)^2$ ($u_\theta(x)$). Binary variable $v_{i,j,j'}$  takes a value of 1 iff $x_{i,j}=x_{i,j'}=1$. Similarly, $r_{j,j'}$ equals 1 iff $x_{1,j}=x_{2,j'}=1$. We formulate the baseline RMP as: 
\begin{align}
\max \quad &u_1+s'\frac{1}{\sqrt{2\pi}} \label{BRMPObj}\\
\text{s.t. } &u_1 = \sum_{j = 1}^n \mu_j  x_{1,j}; \ u_2 = \sum_{j = 1}^n \mu_j  x_{2,j}; \ u_1 \geq u_2 \label{BRMP1}\\
&s = \sum_{i=1}^2\left(\sum_{j=1}^n \sigma^2_j  x_{i,j} + 
2 \sum_{1 \leq j < j' \leq n} \mathrm{cov}(Y_{j} , Y_{j'}) v_{i,j,j'} \right)-2\sum_{j = 1}^{n} \sum_{j' = 1}^n  \mathrm{cov} \left( Y_{j} , Y_{j'}  \right)  r_{j,j'}
\label{BRMP2}\\
&v_{i, j,j'} \leq x_{i,j}; \ v_{i,j,j'} \leq x_{i,j'} \hspace{140pt} \forall j,j'\in \set{1,\dots,n}, \ i \in \set{1,2} \label{BRMP3}\\
&v_{i, j,j'} \geq x_{i,j}+x_{i,j'}-1 \hspace{154pt} 
\forall j,j' \in \set{1,\dots,n}, \ i \in \set{1,2} \label{BRMP4}\\
&r_{j,j'} \leq x_{1,j}; \ r_{j,j'} \leq x_{2,j'} \hspace{149pt} 
\forall j,j' \in \set{1,\dots,n} \label{BRMP5}\\
&r_{j,j'} \geq x_{1,j}+x_{2,j'}-1 \hspace{158pt} 
\forall j,j' \in \set{1,\dots,n} \label{BRMP6}\\
&  \sum_{q=1}^{d}w_q = 1;  \ s' = \sum_{q=1}^{d} \theta_{q+1} w_q \label{Bsrmp:2} \\
& \theta^2_q w_q \leq s \leq \theta^2_{q+1} + \theta^2_{d+1}(1-w_q) 
\hspace{148pt} q = 1, \ldots, d \label{Bsrmp:3} \\
&v \in \set{0,1}^{n \times n \times 2}; \ r \in \set{0,1}^{n \times n}; \ w \in \{0,1\}^{d};  \ x \in \Omega. \label{BRMP7}
\end{align}
The objective function \eqref{BRMPObj} maximizes the upper bounding function defined by Proposition \ref{Pbase}. Constraints \eqref{BRMP1} define the $u$-variables according to equation \eqref{expectedVal} and impose the symmetry breaking condition $u_1 \geq u_2$. Constraint \eqref{BRMP2} imposes $s=\theta(x)^2$ as described by equation \eqref{thetaSquared}. Constraints \eqref{BRMP3}--\eqref{BRMP6} are the McCormick linearization constraints. Constraints~(\ref{Bsrmp:2}) ensure that exactly one interval is chosen for $\theta(x)^2$ and set $s'$ equal to the upper bound of $\theta(x)$ for the interval that $\theta(x)^2$ belongs to. Constraints (\ref{Bsrmp:3}) select the right interval for $\theta(x)^2$. Constraints~\eqref{BRMP7} define the domains of the variables appropriately.

\section{Proofs of Propositions for the Strong Valid Inequalities} \label{app:proofs}

The statements of the propositions are reproduced here for ease of reference.

\setcounter{proposition}{1}
\begin{proposition}
Given~$\delta$ and~$z^{LB}$, let~$\Omega\left(\diff,z^{LB}\right)$ be the set of solutions~$x$ such that~$\diff(x)=\diff$ and~$\mathbb{E} \left[ \max \{ Z_1(x),Z_2(x) \} \right] \geq z^{LB}$. A lower bound~$\underline{\theta}(\delta,z^{LB})$ of~$\theta(x)$ for all $x$ in~$\Omega\left(\diff,z^{LB}\right)$ is given by
\begin{align}
\underline{\theta}(\delta,z^{LB})
= \min_{\theta \geq 0}\left\{ \theta \ | \ \bar{u}+\theta \phi\left(\frac{\delta}{\theta}\right) \geq z^{LB} \right\}.
\end{align}

\end{proposition}

\begin{Proof}
Assume by contradiction that there exists an~$x \in \Omega\left(\diff,z^{LB}\right)$
such that $\theta(x) < \underline{\theta}(\delta,z^{LB})
$. Since $\mathbb{E} \left[ \max \{ Z_1(x),Z_2(x) \} \right] \geq z^{LB}$ we have that 
\begin{equation*}
\mathbb{E}\left[Z_1 (x)\right]  \Phi\left(\frac{\delta}{\theta(x)}\right)+\mathbb{E}\left[Z_2 (x)\right] \Phi\left(\frac{-\delta}{\theta(x)}\right)+\theta(x) \phi\left(\frac{\delta}{\theta(x)}\right) \geq z^{LB}.
\end{equation*}
Because $\bar{u} = \max_{x \in \Omega}\mathbb{E}[Z_1(x)]$, we have
\begin{equation*}
\mathbb{E}\left[Z_1 (x)\right]  \Phi\left(\frac{\delta}{\theta(x)}\right)+\mathbb{E}\left[Z_2 (x)\right]  \Phi\left(\frac{-\delta}{\theta(x)}\right) \leq \bar{u}.
\end{equation*}
As a result,
\begin{equation*}
\bar{u} + \theta(x) \phi\left(\frac{\delta}{\theta(x)}\right) \geq z^{LB},
\end{equation*}
thus contradicting the optimality of~$\underline{\theta}(\delta,z^{LB})$. $\blacksquare$ \end{Proof}

\begin{proposition}
Given~$\diff_1$, $\diff_2$, and~$z^{LB}$ such that $\diff_1 \leq \diff_2$, we have $\underline{\theta}(\delta_1,z^{LB}) \leq \underline{\theta}(\delta_2,z^{LB})$.
\end{proposition}
\begin{Proof}
For $\theta \geq 0$ and fixed~$\delta$, both~$\theta$ and~$\phi\left(\frac{\delta}{\theta}\right)$ are continuous and non-decreasing functions of $\theta$, so we have that $\theta  \phi\left(\frac{\delta}{\theta}\right)$ is also a continuous and non-decreasing function of $\theta$. Therefore, at optimality we have
$
\underline{\theta}\left(\delta_1,z^{LB}\right) \phi\left(\frac{\delta_1}{\underline{\theta}\left(\delta_1,z^{LB}\right)}\right) 
= 
\underline{\theta}(\delta_2,z^{LB}) \phi\left(\frac{\delta_2}{\underline{\theta}\left(\delta_2,z^{LB}\right)}\right) 
= 
z^{LB}-\bar{u}$. 
Assume by contradiction that~$\underline{\theta}(\delta_2,z^{LB}) < \underline{\theta}\left(\delta_1,z^{LB}\right)$; 
if this holds, we must also have  
$
\phi\left(\frac{\delta_2}{\underline{\theta}\left(\delta_2,z^{LB}\right)}\right)
<
\phi\left(\frac{\delta_1}{\underline{\theta}\left(\delta_1,z^{LB}\right)}\right)
$, 
and, consequently,  that 
$ 
\underline{\theta}(\delta_2,z^{LB}) \phi\left(\frac{\delta_2}{\underline{\theta}\left(\delta_2,z^{LB}\right)}\right) 
<  
\underline{\theta}\left(\delta_1,z^{LB}\right) \phi\left(\frac{\delta_1,z^{LB}}{\underline{\theta}\left(\delta_1,z^{LB}\right)}\right)
$, thus contradicting the optimality of~$\underline{\theta}\left(\delta_1,z^{LB}\right)$. 
$\blacksquare$
\end{Proof}

\section{Parameters Used in the Computational Implementation} \label{app:param}
The following are the parameters used in our implementation for the exact algorithm, after fine-tuning with a small subset of instances:
\begin{itemize}
    \item The number of intervals for the discretizations of $\theta(x)$ and $\delta(x)$ are $d=25$ and $l=15$ for the synthetic instances, $d=l=50$ for the machine scheduling application, and $d=50$ and $l=10$ for the DFS application.
    
    \item To obtain upper bounds on $\theta(x)$ and $\delta(x)$ we use the MILP described in Section \ref{sec:Exact} with a time limit of 120 seconds.
    
    \item The large $M$ constant used for constraints \eqref{srmp:7} is the expected makespan of assigning all the jobs to one machine and leaving the other one idle for the scheduling problem. For the other applications we use the maximum profit (score) that can be obtained by a single knapsack (team entry), which we exactly compute via an MILP that solves in less than a second over our testbed.
    
    \item The large $M$ constant used for constraints \eqref{srmp:7.5} is $\theta^2_{d+1} \left(\frac{1}{\sqrt{2\pi}}\right)$. Recall that $\theta^2_{d+1}$ is an upper bound on $\theta(x)$.   
    
    \item The primal heuristic used to obtain a starting feasible solution and lower bound on the objective is as follows. We use the RMP formulation with an additional constraint that requires the model to only select items within the top 90\% expected profit or score. We run this MILP for 60 seconds and recover the best feasible solution found within the time limit.
    
    \item We include the proposed SVIs for solving the synthetic instances and the fantasy football problem. We do not include them in the scheduling application.
    
    \item We stop the algorithm with a tolerance of 0.001, i.e., once $\frac{UB-LB}{UB} < 0.001$ for maximization or $\frac{UB-LB}{LB} < 0.001$ for minimization. 
\end{itemize}

\end{document}